\documentclass[12pt]{article}

\usepackage[latin1]{inputenc}
\usepackage{amsmath,amsthm,amssymb,amscd}

\usepackage{amsmath,amssymb,amsfonts,graphicx,color,fancyhdr}

\newtheorem{thm}{Theorem}[section]
 
 \newtheorem{lem}[thm]{Lemma}
 \newtheorem{prop}[thm]{Proposition}
 \theoremstyle{definition}
 \newtheorem{df}[thm]{Definition}
 \theoremstyle{remark}
 \newtheorem{rem}[thm]{Remark}
 
 \numberwithin{equation}{section}

\def\be#1 {\begin{equation} \label{#1}}
\newcommand{\ee}{\end{equation}}
\def\dem {\noindent {\bf Proof : }}

\def\sqw{\hbox{\rlap{\leavevmode\raise.3ex\hbox{$\sqcap$}}$%
\sqcup$}}
\def\findem{\ifmmode\sqw\else{\ifhmode\unskip\fi\nobreak\hfil
\penalty50\hskip1em\null\nobreak\hfil\sqw
\parfillskip=0pt\finalhyphendemerits=0\endgraf}\fi}

\newcommand{\mb}{\medskip\noindent}
\newcommand{\gb}{\bigskip\noindent}
\newcommand{\R}{\mathbb R}

\newcommand{\N}{\mathbb N}

\newcommand{\C}{\mathbb C}

\newcommand{\B}{\beta}
\newcommand{\T}{r}

\title{Use of abstract Hardy spaces, Real interpolation and Applications to bilinear operators.}
\author{Fr\'ed\'eric Bernicot \\
 \small Universit\'e de Paris-Sud, Orsay et CNRS 8628, 91405 Orsay Cedex, France \\
\small {\em E-mail address:} {Frederic.Bernicot@math.u-psud.fr}}

\begin{document}

\maketitle

\begin{abstract}  

\footnotesize This paper can be considered as the sequel of \cite{BJ}, where the authors have proposed an abstract construction of Hardy spaces $H^1$. They shew an interpolation result for these Hardy spaces with the Lebesgue spaces. Here we describe a more precise result using the real interpolation theory and we clarify the use of Hardy spaces. Then with the help of the bilinear interpolation theory, we then give applications to study bilinear operators on Lebesgue spaces. These ideas permit us to study singular operators with singularities similar to those of bilinear Calder\'on-Zygmund operators in a far more abstract framework as in the euclidean case.

\vspace{0.1in} \footnotesize {\bf Key words}: Hardy spaces, bilinear operators, atomic decomposition, interpolation. 
\vspace{0.1in}

\end{abstract}

\tableofcontents

\section{Introduction.}

The theory of real Hardy spaces started in the 60's, and in the 70's
the atomic Hardy space appeared. Let us recall its definition first
(see \cite{CW}).

\gb Let $(X,d,\mu)$ be a space of homogeneous type and $\epsilon>0$
be a fixed parameter. A function $m\in L^{1}_{loc}(X)$ is called an
$\epsilon$-molecule associated to a ball $Q$ if $\int_X m d\mu=0$ and if
for all $i\geq 0,$
$$ \left( \int_{2^{i+1}Q \setminus 2^i Q} |m|^2 d\mu \right)^{1/2} \leq \mu(2^{i+1}Q)^{-1/2} 2^{-\epsilon i} \textrm{  and  } \left( \int_Q |m|^2 d\mu \right)^{1/2} \leq \mu(Q)^{-1/2}.$$
We call $m$ an atom if in addition we have $ \textrm{supp} (m) \subset Q$. 
Then a function $f$ belongs to $H^1_{CW}(X)$ if there exists a decomposition
$$f =\sum_{i\in\N} \lambda_i m_i  \qquad \mu-a.e, $$
where $m_i$ are $\epsilon$-molecules and $\lambda_i$ are
coefficients which satisfy
$$\sum_{i} |\lambda_i| <\infty.$$
It was proved in \cite{steine} that the whole space $H^1_{CW}(X)$ does not depend on $\epsilon$, as in fact one obtains the same space replacing $\epsilon$-molecules by atoms or $\epsilon'$-molecules with $\epsilon'>0$.

\gb In the Euclidean case ($X=\R^n$ equipped with the Lebesgue measure) this
space has different characterizations, thanks to \cite{steine}~:
\begin{align}
 f\in H^1_{CW}(\R^n)  & \Longleftrightarrow f\in {\mathcal H}^{1} := \left\{f\in L^1(\R^n)\, ;\, \nabla (\sqrt{-\Delta})^{-1}(f) \in L^1(\R^n,\R^n) \right\} \label{Hrond} \\
   & \Longleftrightarrow x \mapsto \sup_{\genfrac{}{}{0pt}{}{y\in\R^n,t>0}{|x-y|\leq t}} \left| e^{-t\sqrt{-\Delta}}(f)(y) \right| \in L^1(\R^n) \label{cha1} \\
  & \Longleftrightarrow x \mapsto \left(\int_{\genfrac{}{}{0pt}{}{y\in\R^n,t>0}{|x-y|\leq t}} \left| t\nabla e^{-t\sqrt{-\Delta}}(f)(y) \right|^2 \frac{dydt}{t^{n+1}} \right)^{1/2} \in L^1(\R^n), \label{cha2}
\end{align}
where $\nabla (\sqrt{-\Delta})^{-1}$ is the Riesz transform. The space ${\mathcal H}^1$ defined by (\ref{Hrond}) was
the original Hardy space of E.M. Stein (see \cite{Se}) and \cite{steine} provided the equivalence with the definition
using the maximal function and the area integral. The link with $H^1_{CW}(\R^n)$ (due to R. Coifman \cite{FS}) comes from the identification of the two dual spaces $(H^1_{CW})^* = ({\mathcal H}^{1})^*$.

\gb The space $H^1_{CW}(X)$ is a good substitute of $L^1(X)$ for many reasons. For instance, Calder\'on-Zygmund operators map $H^1_{CW}(X)$ to $L^1(X)$ whereas they do not map $L^1(X)$ to $L^1(X)$. In addition, $H^1_{CW}(X)$ (and its dual) interpolates with Lebesgue spaces $L^p(X)$, $1<p<\infty$. That is why $H^1_{CW}(X)$ is a good space to extend the scale of Lebesgue spaces $(L^p(X))_{1<p<\infty}$ when $p$ tends to $1$ and its dual $BMO$ when $p$ tends to $\infty$. In addition its atomic decomposition is very useful : for example to check that the set of atoms is sent by a Calder\'on-Zygmund operator in a $L^1$-bounded set, is very easy. That is why we are interested to work with this main property of the atomic (or molecular) structure. We would like also to emphasize that we are more interested by the set of atoms and the vectorial space generated by this collection than by the whole Hardy space.
As we will see, only the behaviour of an operator on the atoms is necessary to use interpolation with Lebesgue spaces. The whole Hardy space is more interesting for example to obtain a characterization of the dual space.

\gb We invite the reader to read \cite{BJ} in order to understand the construction of our Hardy spaces. The work is based on the following remark : there are situations where $H^1_{CW}(X)$ is not the
 right substitute to $L^1(X)$ and there are many works
 where adapted Hardy spaces are defined~: \cite{A1,A,AR,D,DY,D1,D2,DY1,DZ,HM}.
 That is why in \cite{BJ}, we have defined an abstract method to construct Hardy spaces by a molecular (or atomic) decomposition. \\
 In several recent works \cite{GK}, \cite{GK2}, \cite{GT} and \cite{DGY}  X.T. Duong, L. Grafakos, N. Kalton, R. Torres and L. Yan have studied in details some multilinear operators related to multilinear Calder\'on-Zygmund operators on the Euclidean space. \\
 Concerning the linear theory, the abstract Hardy spaces constructed in \cite{BJ} allow us to study linear operators generalizing the study of linear Calder\'on-Zygmund operators. In this paper, we make use of these Hardy spaces to construct a bilinear theory in a most abstract background. Bilinear interpolation theory requires to have a real linear interpolation result. This motives us to use more precisely the ideas of \cite{BJ} to characterize some intermediate spaces between Hardy and Lebesgue spaces for the real interpolation theory (Section \ref{sectionint}).
Then in Section \ref{sectionapp}, we give applications for linear and bilinear operators. By using bilinear interpolation, we will be able to generalize the study of bilinear Calder\'on-Zygmund operators to far more general bilinear operators associated to other cancellations and give examples in Section \ref{sectionexa}.

\section{Definitions and properties of Hardy spaces.}

\gb Let $(X,d,\mu)$ be a space of homogeneous type. We shall write $L^p$ for the Lebesgue space $L^p(X,\R)$ if no confusion arises. Here we are working with real valued functions and we will use "real" duality. We have the same results with complex duality and complex valued functions. \\
By ``space of homogeneous type'' we mean that $d$ is a quasi-distance on the space $X$ and $\mu$ a Borel measure
satisfying the doubling property~: \be{homogene} \exists A>0,  \
\exists \delta>0, \qquad \forall x\in X,\forall r>0, \forall t\geq
1,\qquad \frac{\mu(B(x,tr))}{\mu(B(x,r))} \leq A t^{\delta}, \ee
where $B(x,r)$ is the open ball with center $x\in X$ and radius
$r>0$. We call $\delta$ the homogeneous dimension of $X$. For $Q$ a ball,
and $i\geq 0$, we write $S_i(Q)$ the scaled
corona around the ball $Q$~:
$$ S_i(Q):=\left\{ x,\ 2^{i} \leq 1+\frac{d(x,c(Q))}{r_Q} < 2^{i+1} \right\},$$
where $r_Q$  is the radius of the ball $Q$ and $c(Q)$ its center.
 Note that $S_0(Q)$ corresponds to the ball $Q$ and $S_i(Q) \subset 2^{i+1}Q$ for $i\geq 1$, where $\lambda Q$ is as usual the ball with center $c(Q)$ and radius $\lambda r_Q$.

\mb Before we describe our Hardy spaces, let us recall the definition of Lorentz spaces and give their main properties (for more details see Section 1.4 of \cite{Gra})~:

\begin{df} For  $0< p,s \leq \infty$ we denote $L^{p,s}=L^{p,s}(X)$ the Lorentz space defined by the following norm~:
$$ \|f \|_{L^{p,s}} := \left( \int_0^\infty \left[ t^{1/p} f^*(t) \right]^{s} \frac{dt}{t} \right)^{1/s},$$
where $f^*$ is the decreasing rearrangement of $f$~:
$$f^*(t):= \inf\left\{u>0,\ \mu\left\{|f|>u\right\}\leq t\right\}.$$
\end{df}

\begin{prop} \label{lorentz} 
$1)$ For all exponent $1\leq p\leq \infty$, $L^{p,p}=L^p$.\\
$2)$ For all exponents $0<p,s\leq \infty$, the space $L^{p,s}$ is a metric complete space. Morevover if $p\in(1,\infty]$, $L^{p,s}$ is a Banach space. \\
$3)$ For all exponents $p\in(1,\infty)$ and $ s\in[1,\infty)$, the dual space $(L^{p,s})^*$ is equivalent to the space $L^{p',s'}$.
\end{prop}

\mb We now define the Hardy spaces. Let us denote by $\mathcal{Q}$ the collection of all balls of the space $X$~:
$$ \mathcal{Q}:= \Big\{ B(x,r),\ x\in X, r>0 \Big\}.$$
Let $\B\in(1,\infty]$ be a fixed exponent and $\mathbb{B}:=(B_Q)_{Q\in \mathcal{Q}}$ be a collection of $L^\B$-bounded linear operators, indexed by the collection $\mathcal{Q}$.
We assume that these operators $B_Q$ are uniformly bounded on $L^\B$ : there exists a
constant $0<A'<\infty$
so that~: \be{operh} \forall f\in L^\B ,\ \forall  Q \textrm{ ball}, \qquad \|B_Q(f) \|_{L^\B} \leq A'\|f\|_{L^\B}. \ee

\mb In the rest of the paper, we allow the constants to depend on $A$, $A'$, $\B$ and $\delta$.

\mb We define {\it atoms} and {\it molecules} by using the collection $\mathbb{B}$. We have to think these operators $B_Q$ as the ``oscillation operators'' associated to the ball $Q$.

\begin{df} Let $\epsilon>0$ be a fixed parameter. A function $m\in L^{1}_{loc}$ is called an $\epsilon$-molecule associated to a ball $Q$ if there exists a real function $f_Q$ such that
$$m=B_Q(f_Q),$$
with
$$\forall i\geq 0, \qquad  \|f_Q\|_{L^\B,S_i(Q)} \leq \left(\mu(2^{i}Q)\right)^{-1/\B'} 2^{-\epsilon i}.$$
We call $m=B_Q(f_Q)$ an atom if in addition we have $supp(f_Q) \subset Q$. So an atom is exactly an $\infty$-molecule.
\end{df}

\mb The functions $f_Q$ in this definition are normalized in
$L^1$. It is easy to show that
$$ \|f_Q\|_{L^1} \lesssim 1 \qquad \textrm{and} \qquad \|f_Q\|_{L^\B} \lesssim \mu(Q)^{-1/\B'}.$$
It follows from the $L^\B$-boundedness of the operators $B_Q$ that each molecule belongs to the space $L^\B$. However a molecule is
not (for the moment) in the space $L^1$. Now we are able to define our abstract
Hardy spaces~:

\begin{df} A measurable function $h$ belongs to the molecular
Hardy space $H^1_{\epsilon,mol}$ if there exists a decomposition~:
$$h=\sum_{i\in\N} \lambda_i m_i  \qquad \mu-a.e, $$
where for all $i$, $m_i$ is an $\epsilon$-molecule and $\lambda_{i}$
are real numbers satisfying
$$\sum_{i\in \N} |\lambda_i| <\infty. $$
We equip $H^1_{\epsilon,mol}$ with the norm~:
$$\|h\|_{H^1_{\epsilon,mol}}:= \inf_{h=\sum_{i\in\N} \lambda_i m_i} \sum_{i} |\lambda_i|.$$
Similarly we define the atomic space $H^1_{ato}$ replacing $\epsilon$-molecules by atoms.
\end{df}

\mb In the notations, we forget the exponent $\B$ and the collection ${\mathbb B}$. As we will explain in Subsection \ref{linear}, we need to use a smaller space~:

\begin{df} \label{hardyF} According to the collection ${\mathbb B}$, we introduce the set $H^1_{F,\epsilon, mol}\subset H^1_{\epsilon, mol} \cap L^\B$, given by the finite sum of $\epsilon$-molecules with the following norm
$$ \|f\|_{H^1_{F,\epsilon, mol}} := \inf_{f=\sum_i \lambda_i m_i} \sum_i |\lambda_i|. $$
We take the infimum over all the finite molecular decompositions.
Similarly we define the atomic space $H^1_{F,ato}$.
\end{df}

\mb As we will see in Subsection \ref{linear}, it is quite easy to estimate the behaviour of an operator on the whole collection of atoms or molecules. Then by linearity, we can control the operator on the previous ``finite'' Hardy spaces with the corresponding norm. However to extend the operator on the whole Hardy space (in a continuous way with its smaller norm) is an abstract problem which seems quite difficult and requires some extra assumptions (see Subsection \ref{linear}). Fortunately, we will see that we do not need to study if the operator can be extended or not, its behaviour on the sets of atoms will be sufficient.

\mb To understand in a better way our definition, we refer the reader to Section 3 of \cite{BJ}, where we compare our Hardy spaces with some already studied Hardy spaces. Let us make some remarks.

\begin{rem}
\mb $1-)$ For the Hardy space $H^1_{\epsilon, mol}$, we only ask that the decomposition
$$h(x) =\sum_{i\in\N} \lambda_i m_i(x) $$
is well defined for almost every $x\in X$. So the assumption is very
weak and it is possible that the measurable function $h$ does not
belong to $L^1_{loc}$. It is not clear whether these abstract normed vector
spaces are complete. The problem is that we do not know whether
the decompositions of $h$ converge absolutely. \\
$2-)$ We have the following continuous embeddings~:
\be{inclusio} \forall\, 0<\epsilon<\epsilon',\  \qquad  H^1_{F,ato} \hookrightarrow H^1_{ato} \hookrightarrow H^1_{\epsilon',mol} \hookrightarrow H^1_{\epsilon,mol}. \ee
In fact the space $H^1_{ato}$ corresponds to the space $H^1_{\infty,mol}$.
 For $0<\epsilon<\epsilon'<\infty$ the space $H^1_{\epsilon',mol}$ is
 dense in $H^1_{\epsilon, mol}$. In the general case,
 it seems to be very difficult to study the dependence of $H^1_{\epsilon, mol}$ with the parameter $\epsilon$ and we will not study this question here. Similarly the dependence of the Hardy spaces on the exponent $\B$ is an interesting question, but seems difficult and will not be studied here. \\
$3-)$ The norm on the molecular spaces $H^1_{\epsilon, mol}$ and on the finite molecular spaces $H^1_{F,\epsilon, mol}$ may not be equivalent (see Subsection \ref{linear} and a counterexample of Y. Meyer for the Coifman-Weiss space in \cite{Meyer}).
\end{rem}

\begin{rem} \label{lemdense}
 We have seen that each molecule is an $L^\B$ function. So it is obvious that $H^1_{F,\epsilon, mol} \subset L^\B \cap H^1_{\epsilon, mol}$ is dense in $H^1_{\epsilon, mol}$ and $H^1_{F,ato} \subset L^\B \cap H^1_{ato}$ is dense in $H^1_{ato}$.
\end{rem}

\mb To work with a vector normed space, we often need the completness of this one. The following proposition gives us some conditions to get the completness of our Hardy spaces.

\begin{prop} \label{banach} Take $\epsilon \in ]0,\infty]$ and assume that the space $H^1_{\epsilon, mol}$ is continuously embedded in $L^1_{loc}$. Then $H^1_{\epsilon, mol}$ is a Banach space.
 \end{prop}

\dem We have just to verify the
 completeness. The proof is easy by
using the following well-known criterion : for $\epsilon>0$,
$H^1_{\epsilon,mol}$ is a Banach space if for all sequences
$(h_i)_{i\in\N}$ of $H^1_{\epsilon,mol}$ satisfying
$$\sum_{i\geq 0} \left\| h_i \right\|_{H^1_{\epsilon,mol}} <\infty, $$
the series $\sum h_i$ converges in the Hardy space
$H^1_{\epsilon,mol}$. This is true because each molecular
decomposition is absolutely convergent in $L^1_{loc}$-sense. We therefore
define the series $\sum_{i} h_i$ as a measurable function in $L^1_{loc}$.
Finally, it is easy to prove the convergence of the series for the $H^1_{\epsilon,mol}$
norm. \findem

\mb The next proposition explains that under some ``fast'' decays for the operators $B_Q$, the Hardy spaces are included in $L^1$.

\begin{prop} Assume that the operators $B_Q$ satisfy that for $M>\delta/\B'$ a large enough integer, there exists a constant $C$ such that for all $i,k\geq 0$
\be{decay}  \forall f\in L^\B,\ \textrm{supp}(f) \subset 2^i Q \qquad \left\| B_Q(f) \right\|_{L^\B,S_k(2^{i}Q)} \leq C 2^{-Mk}
\|f\|_{L^\B,2^{i}Q}. \ee
Then the following inclusions hold~:
$$ \forall \epsilon>0, \qquad H^{1}_{ato} \hookrightarrow H^{1}_{\epsilon, mol} \hookrightarrow L^1.$$
Consequently, the Hardy spaces are Banach spaces.
\end{prop}

\dem We claim that all $\epsilon$-molecules (and atoms) are bounded in $L^1$. In fact, using (\ref{decay})
\begin{align*}
\left\| B_Q(f_Q) \right\|_{L^1} & \leq  \sum_{i\geq 0} \left\| B_Q(f_Q{\bf 1}_{S_i(Q)}) \right\|_{L^1} \leq \sum_{i\geq 0} \sum_{k\geq 0} \left\| B_Q(f_Q{\bf 1}_{S_i(Q)}) \right\|_{L^1,S_k(2^{i}Q)} \\
 & \lesssim \sum_{i\geq 0} \sum_{k\geq 0} \mu(2^{i+k} Q)^{1/\B'} \left\| B_Q(f_Q{\bf 1}_{S_i(Q)}) \right\|_{L^\B,S_k(2^{i}Q)} \\
 & \lesssim \sum_{i\geq 0} \sum_{k\geq 0} \mu(2^{i+k} Q)^{1/\B'} 2^{-kM} \left\|f_Q \right\|_{L^\B,S_i(Q)} \\
 & \lesssim \sum_{i\geq 0} \sum_{k\geq 0} \mu(2^{i+k} Q)^{1/\B'} 2^{-Mk} \mu(2^iQ)^{-1/\B'} 2^{-\epsilon i} \\
 & \lesssim \sum_{i\geq 0} \sum_{k\geq 0} 2^{k\delta/\B'} 2^{-Mk} 2^{-\epsilon i} \lesssim 1.
\end{align*}
 Here we use the estimates for $f_Q$, the doubling property of $\mu$ and the fact that $M$ is large enough ($M>\delta/\B'$ works with $\B'<\infty$). Thus we obtain that all $\epsilon$-molecules are bounded in $L^1$, and we can deduce the embedding from the definition of the Hardy spaces. \findem

\mb We have seen that all the molecular spaces contain the atomic space. In the next section, we will study the intermediate spaces, obtained by real interpolation between the Hardy spaces and the Lebesgue spaces. We will see that under reasonnable assumptions, they ``do not depend'' on the considered Hardy space.

\gb
\section{Real Interpolation between Hardy and Lebesgue spaces.}

\label{sectionint}

\subsection{Preliminaries about Real interpolation theory.}

Let us begin to remember the definition of the real interpolation theory (see for more details the book \cite{BL}).

\begin{df} Let $E$ and $F$ be two vector normed spaces. For every $a\in E+F$ and every $t>0$, we define the $K$-functional as
 $$ K(t,a,E,F):= \inf_{\genfrac{}{}{0pt}{}{a=e+f}{ \genfrac{}{}{0pt}{}{e\in E}{f\in F}}} \|e\|_E + t\|f\|_F.$$
For every $a\in E \cap F$ and every $t>0$, we define the $J$-functional as
 $$ J(t,a,E,F):= \max\big\{ \|a\|_E,\ t\|f\|_F \big\}.$$
In addition we say that the couple $(E,F)$ is compatible if the space $E\cap F$ is dense in $E+F$.
\end{df}

\mb Now with these two functionals, we can define some particular intermediate spaces.

\begin{df} \label{definter} Let $E$ and $F$ be two vector normed spaces. Then for $\theta\in (0,1)$ and $s\in [1,\infty]$, we denote by $(E,F)_{\theta,s,K}$ and $(E,F)_{\theta,s,J}$ the following spaces
 $$ (E,F)_{\theta,s,K} := \left\{ a\in E+F,\ \|a\|_{(E,F)_{\theta,p,K}}:= \left(\int_0^\infty \left(t^{-\theta} K(t,a,E,F)\right)^s \frac{dt}{t} \right)^{1/s} <\infty \right\} $$
and similarly
 $$ (E,F)_{\theta,s,J} := \left\{ a\in E+F,\ \|a\|_{(E,F)_{\theta,s,J}} <\infty \right\}, $$
with 
$$ \|a\|_{(E,F)_{\theta,s,J}} := \inf_{a=\int_0^\infty u(t) \frac{dt}{t}} \left(\int_0^\infty \left(t^{-\theta} K(t,u(t),E,F)\right)^s \frac{dt}{t} \right)^{1/s}. $$
\end{df}

\mb We recall the well-known results about these spaces (for the proofs of the following results and any details about the real interpolation, we refer to the book \cite{BL} of J. Bergh and J. L\"ofstr\"om).

\begin{thm}[Equivalence Theorem] Let $E$ and $F$ be two vector normed spaces. For all $\theta\in(0,1)$ and $s\in[1,\infty]$, the two spaces $(E,F)_{\theta,s,K}$ and $(E,F)_{\theta,s,J}$ are equal with equivalent norms. From now, we denote these spaces $(E,F)_{\theta,s}$. In addition this space is an intermediate space, that is
 $$ E\cap F \hookrightarrow (E,F)_{\theta,s} \hookrightarrow E+F $$
with continuous embeddings.
\end{thm}

\mb We will use the link between the spaces and their completed spaces~:

\begin{thm} \label{complet} Let $E$ and $F$ be two vector normed spaces. We set $\overline{E}$ and $\overline{F}$ for the completed spaces. Then for all $\theta\in(0,1)$ and $s\in[1,\infty]$, the two spaces $\overline{(E,F)_{\theta,s}}$ and $(\overline{E},\overline{F})_{\theta,s}$ are equal (with equivalent norms). So we deduce that 
$$ (E,F)_{\theta,s} =\left\{u\in E+F,\ \|u\|_{(\overline{E},\overline{F})_{\theta,s}}<\infty \right\}.$$
\end{thm}

\mb Let us recall the following notion due to \cite{KP}~:

\begin{df} The couple $(E,F)$ is said to be a doolittle couple if the ``diagonal'' space $E\cap F\simeq \{(x,x),x\in E\cap F \}$ is closed in $E\times F$. 
\end{df}

\mb For the interpolation of dual spaces, we have the following theorem~:

\begin{thm}[Duality Theorem] \label{thmdual} Assume that the couple $(E,F)$ is a compatible (or a doolittle) couple of Banach spaces then for all $\theta\in(0,1)$ and $s\in[1,\infty]$ we have
 $$ \left[(E,F)_{\theta,s} \right]^* = (E^*,F^*)_{\theta,s'}, $$ 
where 
$$\frac{1}{s}+\frac{1}{s'}=1.$$
\end{thm}

\mb The proof is given in \cite{BL} for a compatible couple and in \cite{KP} for a doolittle couple. The completeness is important to use Hahn-Banach Theorem in order to obtain this equivalence.

\subsection{Interpolation between Hardy and Lebesgue spaces.}

\mb After recalling these results, we want to study the real interpolation between our Hardy spaces and Lebesgue spaces. Let $\B\in(1,\infty]$ be always fixed, we work with the Hardy space $H^1$ equal to one of the following Hardy spaces : $H^1_{F,ato}$, $H^1_{ato}$, $H^1_{F,\epsilon, mol}$ or $H^1_{\epsilon, mol}$. The space $H^1$ may be the completed space of $H^1_{F,ato}$ or $H^1_{F,\epsilon, mol}$ too. We will need the following definitions~:

\begin{df} We set $A_Q$ for the operator $Id-B_Q$. For $\sigma\in [1,\infty]$ we define the maximal operator~:
\be{opeM} \forall x\in X, \qquad M_{\sigma}(f)(x):= \sup_{\genfrac{}{}{0pt}{}{Q \textrm{ball}}{x\in Q}}\  \left( \frac{1}{\mu(Q)} \int_Q \left|A_Q^*(f)\right|^\sigma d\mu  \right) ^{1/\sigma}, \ee
where $A_Q^*$ is the adjoint operator. \\
For all $s>0$, we define a maximal sharp function adapted to our operators~:
$$ \forall x\in X, \qquad M^\sharp_s (f)(x):= \sup_{\genfrac{}{}{0pt}{}{Q \textrm{ball}}{x\in Q}} \left(\frac{1}{\mu(Q)} \int_{Q} \left|B_Q^* (f)(z) \right|^{s} d\mu(z) \right)^{1/s}.$$
The standard maximal ``Hardy-Littlewood'' operator is defined for $s>0$ by
$$\forall x\in X, \qquad M_{HL,s}(f)(x):=\sup_{\genfrac{}{}{0pt}{}{Q \textrm{ball}}{x\in Q}} \left(\frac{1}{\mu(Q)} \int_{Q} \left|f(z) \right|^{s} d\mu(z)
\right)^{1/s}.$$
\end{df}

\begin{lem} \label{lem1}  For each ball $Q$ of $X$, the operator $B_Q^*$ which is defined on $L^{\B'}$ can be extended to an operator acting on $(H^1)^*$ to ${L^{1}}_{loc}$. Keeping the same notation $B_Q^*$ for the extension, we have that for all $s \in[1,\B']$
$$\forall f\in (H^1)^*, \qquad \left\| M^\sharp_s(f) \right\|_\infty \leq \|f\|_{(H^1)^*}.$$
\end{lem}

\dem Let us fix an element $f\in (H^1)^*$ and a ball
$Q$. For all function $h$, supported
in $Q$ and normalized by $\|h\|_{\B}=1$ then we set $\phi_Q:=
\mu(Q)^{-1/\B'}h$. Then it is obvious that $m=B_Q(\phi_Q)$ is an
atom. By duality, we obtain
$$\left| \langle f, B_Q(\phi_Q) \rangle \right| \leq \|f\|_{(H^1)^*}.$$
Therefore
$$ \forall h\in L^{\B}(Q),\ \|h\|_{L^{\B}}=1, \qquad \left| \langle f, B_Q(h) \rangle \right| \leq \mu(Q)^{1/\B'} \|f\|_{(H^1)^*}\|h\|_{L^\B(Q)}.$$
Hence we can define, $B_Q^*(f)\in L^{\B'}(Q)=(L^{\B}(Q))^*$ such that
$$ \langle B_Q^*(f), h \rangle := \langle f, B_Q(h) \rangle.$$ 
We obtain also the estimate
$$\left\| B_Q^*f \right\|_{L^{\B'}(Q)} \leq  \mu(Q)^{1/\B'} \|f\|_{(H^1)^*},$$
which concludes the proof. \findem

\mb The dual space $(H^1)^*$ is always a (Banach) normed vector space, and we can use real interpolation with an $L^q$ space. For the $K$-functional, we prove ~:

\begin{prop} \label{thm:bonlambda} Assume that for $\sigma> \B'$, the operator $M_\sigma$ is bounded by $M_{HL,\B'}$.
Then for $q\in[\B',\infty)$ there exist a constant $c_1=c_1(q)$ and for all $\kappa\geq 2$ an other constant $c_2=c_2(q,\kappa)\geq 1$ such that for all $t>0$, $s\in [1,\B']$, we have the following estimates for the $K$ functional~: for every function $f\in L^{\B'}$
\begin{align} 
 K(t,f,L^q,(H^1)^*) \geq c_1 K(t,{M_s}^\sharp(f),L^{q,\infty},L^\infty),  \label{ineq1} 
\end{align}
and
\begin{align} 
\lefteqn{ K(t,M_{HL,q}(f),L^{q,\infty},L^\infty) \leq c_2 \Big[K(t,f,L^q,(H^1)^*)  } & & \nonumber \\
 & &   + \kappa^{1-\sigma/q}  K(\kappa^{\sigma/q} t,M_{HL,q}(f),L^{q,\infty},L^\infty)+ t{\bf 1 }_{t^q \lesssim \mu(X)} \|f\|_{L^{\B'}} {\bf 1}_{\mu(X)<\infty} \Big]. \label{ineqK2}
\end{align}
If $\mu(X)<\infty$, we allow the constant $c_2$ to depend on the measure $\mu(X)$.
\end{prop}

\mb The proof of this proposition requires the following lemma.

\begin{lem} \label{good}  Suppose that $M_\sigma$ is bounded by $M_{HL,\B'}$  with $\sigma> \B'$. Then for all $\rho$ large enough and for all $\gamma<1$ we have a ``good lambdas inequality''~: for all $\lambda>0$ (or $\lambda \gtrsim \|h\|_{L^{\B'}}$ if $\mu(X)<\infty$) 
$$ \mu\left( \left\{ M_{HL,\B'}(h) >\rho\lambda, M_{\B'}^\sharp(h)(x)\leq \gamma\lambda \right\} \right) \lesssim \left( \gamma^{\B'} \rho^{-\B'}+ \rho^{-\sigma} \right) \mu\left( \left\{ M_{HL,\B'}(h)>\lambda \right\} \right).$$
The implicit constant does not depend on $\rho$ and $\gamma$.
 In addition, for all exponent $q\in (\B',\sigma)$ there exists a constant $c=c(q)$ such that for all function $h\in L^{\B'}$
\be{equiv}  c\|h\|_{L^q} \leq \|M_{\B'}^\sharp(h) \|_{L^q} + {\bf 1}_{\mu(X)<\infty}\|h\|_{L^{\B'}}  \leq c^{-1} \|h\|_{L^q}. \ee
\end{lem}

\dem This lemma is a consequence of Theorem 3.1 in \cite{AM}. With its notations, take $F=|h|^{\B'}$ and for all balls $Q$
$$ G_Q= c_{\B'}|B_Q^*h|^{\B'} \textrm{    and    } H_Q= c_{\B'}|A_Q^*h|^{\B'},$$
where $c_{\B'}$ is a constant such that $(u+v)^{\B'} \leq c_{\B'}(u^{\B'}+v^{\B'})$.
Then for all $x\in Q$, we get $F(x)\leq G_Q(x)+H_Q(x)$. For all balls $Q \ni x$, we have
$$ \frac{1}{\mu(Q)} \int_Q G_Q d\mu \leq c_{\B'} M_{\B'}^\sharp(h)(x)^{\B'}:=G(x).$$
By assumption, for all balls $Q \ni x$ and $\overline{x}\in Q$, we obtain
\begin{align*} \left( \frac{1}{\mu(Q)} \int_{Q} |H_Q|^{\sigma/\B'} d\mu \right)^{\B'/\sigma} & \leq c_{\B'}^{\B'/\sigma} \left( \frac{1}{\mu(Q)} \int_Q |A_Q^*(h)|^\sigma d\mu \right)^{\B'/\sigma} \leq c_{\B'}^{\B'/\sigma} M_{\sigma}(h)(\overline{x})^{\B'} \\
& \lesssim  M_{HL,\B'}(h)(\overline{x})^{\B'} =  M_{HL,1}(F)(\overline{x}).
\end{align*}
Applying Theorem 3.1 of \cite{AM} (proved for spaces of homogeneous type in Section 5 of \cite{AM}), we obtain the desired ``good lambdas inequality''~:
$$ \mu\left( \left\{ M_{HL,\B'}(h)^{\B'} >\rho\lambda, M_{\B'}^\sharp(h)(x)^{\B'}\leq \gamma\lambda \right\} \right) \lesssim \left(\gamma \rho^{-1} + \rho^{-\sigma/\B'}\right) \mu\left( \left\{ M_{HL,\B'}(h)^{\B'}>\lambda \right\} \right).$$
By the same theorem, we also have that for all exponent $t\in[1,\sigma/\B')$ and for $F\in L^1$ (which is equivalent to $h\in L^{\B'}$)
\be{cont11} \|F\|_{L^t} \leq \| M_{HL,1}(F) \|_{L^t} \lesssim \| M_{\B'}^\sharp(h)^{\B'} \|_{L^t} = \|M_{\B'}^\sharp(h) \|_{L^{t\B'}}^{\B'}. \ee
Therefore for all $q\in[\B',\sigma)$
$$ \| f \|_{L^q} \lesssim \|M_{\B'}^\sharp(h) \|_{L^q}.$$
The extra term, which appears when $\mu(X)<\infty$, is explained in Section 5 of \cite{AM}.
The other inequality for (\ref{equiv}) is much more easy and is a direct consequence of the fact that $M_{HL,\B'}$ is $L^q$-bounded and
$$ M_{\B'}^\sharp(h) \lesssim M_{HL,\B'}(h)+M_{\B'}(h) \lesssim M_{HL,\B'}.$$
\findem

\begin{rem} Assume that $\mu(X)=\infty$. Then for $\B'< q<\sigma$, we obtained a ``Fefferman-Stein'' inequality~:
$$  c\|h\|_{L^q} \leq \|M_{\B'}^\sharp(h) \|_{L^q} \leq c^{-1} \|h\|_{L^q}. $$
In fact the proof shows that the right inequality is true for all $q>\B'$.
\end{rem}

\gb
Now we can prove Proposition \ref{thm:bonlambda}

\mb
\dem By definition of the $K$-functional, we have~:
 $$ K(t,f,L^q,(H^1)^*):= \inf_{\genfrac{}{}{0pt}{}{f=g+h}{ \genfrac{}{}{0pt}{}{g\in L^q}{h\in (H^1)^*}}} \|g\|_{L^q} + t\|h\|_{(H^1)^*}.$$
We want to use the maximal function $M_s^\sharp$. From the assumption on the maximal function $M_\sigma$, it is easy to see that for every function $\phi$ 
$$ M^\sharp_s(\phi) \lesssim M_{HL,s}(\phi) + M_s(\phi) \lesssim M_{HL,q}(\phi).$$
So $M^\sharp_s$ is of weak type $(q,q)$ as $s\leq \B'\leq q$. We have also~:
\begin{align*} 
 K(t,f,L^q,(H^1)^*) & \gtrsim \inf_{\genfrac{}{}{0pt}{}{f=g+h}{ \genfrac{}{}{0pt}{}{g\in L^q}{h\in (H^1)^*}}} \|M_s^\sharp(g)\|_{L^{q,\infty}} + t\|h\|_{(H^1)^*} \\
 & \gtrsim \inf_{\genfrac{}{}{0pt}{}{f=g+h}{ \genfrac{}{}{0pt}{}{M_s^\sharp(g)\in L^{q,\infty}} {h\in (H^1)^*} } }   \|M_s^\sharp(g)\|_{L^{q,\infty}}+t\|h\|_{(H^1)^*}.
\end{align*}
By the same way, using Lemma \ref{lem1}, we have
\begin{align*} 
 K(t,f,L^q,(H^1)^*) & \gtrsim \inf_{\genfrac{}{}{0pt}{}{f=g+h}{ \genfrac{}{}{0pt}{}{M_s^\sharp(g)\in L^{q,\infty}}{h\in (H^1)^*}}} \|M_s^\sharp(g)\|_{L^{q,\infty}} + t\|M_s^\sharp(h)\|_{L^{\infty}}.
\end{align*}
Now we note that for $f=g+h$, we have $M_s^\sharp(f) \leq M_s^\sharp(g) + M_s^\sharp(h)$ and so we can conclude 
$$ M_s^\sharp(h) \geq \left|M_s^\sharp(f) - M_s^\sharp(g)\right|.$$
We get also
\begin{align} 
 K(t,f,L^q,(H^1)^*) & \gtrsim \inf_{\genfrac{}{}{0pt}{}{f=g+h}{ \genfrac{}{}{0pt}{}{{M_s}^\sharp(g)\in L^{q,\infty}}{h\in (H^1)^*}}} \|{M_s}^\sharp(g)\|_{L^{q,\infty}} + t\| {M_s}^\sharp(f) - {M_s}^\sharp(g) \|_{L^\infty}  \nonumber \\
& \gtrsim \inf_{\genfrac{}{}{0pt}{}{{M_s}^\sharp(g)\in L^{q,\infty}}{{M_s}^\sharp(f) - {M_s}^\sharp(g) \in L^\infty}} \|{M_s}^\sharp(g)\|_{L^{q,\infty}} +t\| {M_s}^\sharp(f) - {M_s}^\sharp(g) \|_{L^\infty}  \nonumber \\
& \gtrsim K(t,{M_s}^\sharp(f),L^{q,\infty},L^\infty). \label{relK} 
\end{align}
So we have proved the first desired inequality (\ref{ineq1}). Now we will use the fact that we exactly know the $K$-functional for the Lebesgue spaces. Let us remember the following result (see \cite{BL}, p109)~:
\be{equiK} K(t,\phi,L^{q,\infty},L^{\infty}) \simeq \left[\sup_{0<u} u \left| \left\{ v\in[0,t^q],\ \phi^*(v)>u \right\} \right|^{1/q}\right]:=\left\|\phi^* \right\|_{L^{q,\infty}([0,t^q])}, \ee
where we write $\phi^*$ for the decreasing rearrangement function
$$ \phi^*(v):= \inf \{ \lambda >0,\ \mu(|\phi|>\lambda) \leq v \}.$$ 
We have also to estimate the function $\left[M_{\B'}^\sharp(f)\right]^*$. By Lemma \ref{good}, we obtain~: for all $\lambda>0$ (or $\lambda \gtrsim \|f\|_{L^{\B'}}$ if $\mu(X)<\infty$)
$$ \mu\left( \left\{ M_{HL,q}(f)>\rho\lambda,\ M_{\B'}^\sharp(f)\leq  \gamma \lambda \right\} \right) \lesssim \left(\gamma^{\B'} \rho^{-\B'}+\rho^{-\sigma}\right) \mu\left( \left\{ M_{HL,q}(f)>\lambda\right\}\right).$$
Let $\delta:=\gamma^{\B'} \rho^{-\B'}+\rho^{-\sigma}$. We deduce that (with $C\geq 1$ a constant independent on the important parameters)
$$ \mu\left( \left\{ M_{HL,q}(f)>\rho\lambda \right\} \right) \leq C \mu\left( \left\{ M_{\B'}^\sharp(f)> \gamma \lambda \right\} \right)  + C\delta \mu\left( \left\{ M_{HL,q}(f)>\lambda\right\}\right).$$
Then it follows easily that
$$\left[ M_{HL,q}(f)\right]^*(s) \leq \rho\max\left\{ \gamma^{-1} \left[M_{\B'}^\sharp(f)\right] ^*(s/C), \left[M_{HL,q}(f)\right]^*(\delta^{-1}s/C) \right\} + \lambda_0{\bf 1}_{\mu(X)<\infty},$$ 
where $\lambda_0$ satisfies $\lambda_0 \lesssim \|f\|_{L^{\B'}}$. To prove (\ref{ineqK2}), we distinguish 2 cases~: \\
$1-)$ First case : assume that $\mu(X)=\infty$. \\
With (\ref{relK}) and (\ref{equiK}), we obtain
\begin{align} 
 \lefteqn{K(t,M_{HL,q}(f),L^{q,\infty},L^\infty) \simeq \left\|M_{HL,q}(f)^* \right\|_{L^{q,\infty}([0,t^q])}} & & \nonumber  \\ 
 & & \lesssim  \rho\gamma^{-1} \left\| M_{\B'}^\sharp(f) ^*(\cdot/C) \right\|_{L^{q,\infty}([0,t^q])} + \rho\left\| M_{HL,q}(f) ^*(\delta^{-1}\cdot/C) \right\|_{L^{q,\infty}([0,t^q])} \label{autre} \\
 & & \lesssim  \rho\gamma^{-1} \left\| M_{\B'}^\sharp(f) ^*(\cdot) \right\|_{L^{q,\infty}([0,t^q])} + \rho\delta^{1/q} \left\| M_{HL,q}(f) ^* \right\|_{L^{q,\infty}([0,\delta^{-1}t^q])}. \nonumber
\end{align}
In all these last estimates, the implicit constant depends on the constant $C\geq 1$ and is uniform on the two important constants $\rho\geq 2$ and $\gamma\leq 1$. 
By using the equivalence (\ref{equiK}) and (\ref{relK}), we obtain the desired inequality~:
\begin{align*} 
 K(t,M_{HL,q}(f),L^{q,\infty},L^\infty) \lesssim  \rho\gamma^{-1} K(t,f,L^q,(H^1)^*)
 + \rho\delta^{1/q} K(\delta^{-1/q} t,M_{HL,q}(f),L^{q,\infty},L^\infty).
\end{align*}
We choose $\gamma=\rho^{-\sigma/\B'+1}$ small enough to have $\delta\simeq 2\rho^{-\sigma}$ and then we take $\kappa \simeq \rho$. \\
$2-)$ Second case : $\mu(X)<\infty$. \\
In this case, there is an extra term in (\ref{autre}) which corresponds to 
$$  \|\lambda_0 \|_{L^{q,\infty}([0,t^q])} .$$
This term is bounded by $t\lambda_0$ and one has just to consider $t^q \leq \mu(X)$. Using $ \lambda_0 \lesssim \|f\|_{L^{\B'}}$, we obtain the result with an implicit constant which depends on $\mu(X)$. The previous study holds with the same arguments for the main terms. The proof is also completed.
\findem 

\gb Having good estimates for the interpolation $K$ functional between the Hardy space $H^1$ and the Lebesgue spaces, we now prove the following result ~:

\begin{thm} \label{thm:interpolation} Assume that for $\sigma> \B'$, the operator $M_\sigma$ is bounded by $M_{HL,\B'}$.
Let $\theta$ be in $(0,1]$, $s\in[1,\infty]$, $q\in[\B',\infty)$ and $t$ satisfying
$$ \frac{1-\theta}{q} = \frac{1}{t}>\frac{1}{\sigma}.$$ 
Then there exists a constant $c=c(\theta,s,\sigma,q)$ such that for all function $f\in L^{\B'}\cup L^q$~: if $f\in (L^q,(H^1)^*)_{\theta,s}$ then $f\in L^{t,s}$ and
$$ \|f\|_{L^{t,s}} \leq c \left[\|f\|_{(L^q,(H^1)^*)_{\theta,s}} + \|f\|_{L^{r'}}{\bf 1}_{\mu(X)<\infty}\right].$$
\end{thm}

\dem To compute the norm of the intermediate space, we have to integrate the $K$-functional as decribed in Definition \ref{definter}~:
$$ \|f\|_{(L^q,(H^1)^*)_{\theta,s}} = \left(\int_0^\infty \left(t^{-\theta} K(t,f,L^q,(H^1)^*)\right)^s \frac{dt}{t} \right)^{1/s}. $$
We will only deal with the interesting case $\mu(X)=\infty$. The other case uses the same arguments with an extra term which is very easy to control (as $t$ is bounded).
We use the inequality obtained in Theorem \ref{thm:bonlambda}.
$$ K(t,M_{HL,q}(f),L^{q,\infty},L^\infty) \lesssim  K(t,f,L^q,(H^1)^*) + \kappa^{1-\sigma/q}  K(\kappa^{\sigma/q} t,M_{HL,q}(f),L^{q,\infty},L^\infty). $$
By a formally integration with a change of variable, we get
\be{formaint} \|M_{HL,q}(f)\|_{(L^{q,\infty},L^\infty)_{\theta,s}} \lesssim \|f\|_{(L^q,(H^1)^*)_{\theta,s}} + \kappa^{1-\sigma/q+\theta\sigma/q} \|M_{HL,q}(f)\|_{(L^{q,\infty},L^\infty)_{\theta,s}}. \ee
So for $\kappa$ large enough as $1-\sigma/q+\theta\sigma/q<0$, we obtain that
$$ \|M_{HL,q}(f)\|_{(L^{q,\infty},L^\infty)_{\theta,s}} \lesssim \|f\|_{(L^q,(H^1)^*)_{\theta,s}}.$$
Here we have assumed that the left side of (\ref{formaint}) is finite to conclude.
In fact, by the same arguments as used in Theorem 3.1 of \cite{AM}, it suffices to have that $f\in L^{\B'}$ or $f\in L^q$ to prove this last inequality. In addition the arguments in Section 5 of \cite{AM} allow us to obtain the extra term when $\mu(X)<\infty$. The proof is also finished since we know (see \cite{BL}) that 
$$ (L^{q,\infty},L^\infty)_{\theta,s}= L^{t,s}. $$
\findem 

\mb Now by duality, we look for a real interpolation result for our Hardy space.

\begin{prop} \label{thm:inter} Assume that the Hardy space $H^1$ is complete and satisfies ~: 
$$ H^1 \hookrightarrow L^1.$$
Assume that $\mu(X)=\infty$ and that for $\sigma> \B'$, the operator $M_\sigma$ is bounded by $M_{HL,\B'}$.
Then for all $\theta\in (0,1),s\in[1,\infty)$, for all exponents $p\in(1,\B)$ and $t\in(1,\infty)$ such that 
$$ \frac{1}{t}= \frac{1-\theta}{p} + \theta <\frac{1}{\sigma'},$$
we have the following equivalence : 
\be{injection} L^{t,s} \simeq (L^{p},H^1)_{\theta,s}.\ee
\end{prop}

\dem First we note that under our assumption, the Hardy space $H^1$ is a Banach space (due to Proposition \ref{banach}) if $H^1=H^1_{ato}$ or $H^1_{\epsilon, mol}$. Now the result is well known for $H^1=L^1$ so we already have the inclusion 
\be{inclu} E:=(L^{p},H^1)_{\theta,s} \hookrightarrow L^{t,s} . \ee
Let us write $E:=(L^{p},H^1)_{\theta,s}$ and $(F,\|\ \|_{L^{t,s}})$ be the closure of $E$ for the $L^{t,s}$ norm. So $F$ is a closed subspace of $L^{t,s}$. Assume that $F \subsetneq L^{t,s}$. Then by Hahn-Banach Theorem (and the duality result in Proposition \ref{lorentz}), there exists a function $\phi \in L^{t',s'}$ with $\|\phi\|_{L^{t',s'}}=1$ such that $\phi=0$ in $F^*$.
We have $\phi \in F^*$ and so by (\ref{inclu}), we have that $\phi \in E^*$. \\
In addition, we claim that $(L^p,H^1)$ is a doolittle couple of Banach spaces. We have to check that the (diagonal) space $L^p \cap H^1$ is a closed sub-space of $H^1 \times L^p$. Take a sequence $(x_n,x_n)$ which converges to $(x,y)$ in $H^1 \times L^p$. From the above assumption, $x_n$ converges to $x$ in $L^1$-sense and to $y$ in $L^p$ sense. We deduce that $x=y\in H^1 \cap L^p$. Using the duality Theorem (Theorem \ref{thmdual}), we obtain that
$$ E^* = (L^{p'},(H^1)^*)_{\theta,s'}.$$
Theorem \ref{thm:interpolation} gives us also that 
\be{propdual} 1=\|\phi\|_{L^{t',s'}} \lesssim \|\phi\|_{E^*}. \ee
By the embedding (\ref{inclu}) and the duality properties of Banach spaces, we have
$$ \|\phi\|_{E^*} = \sup_{\genfrac{}{}{0pt}{}{e\in E}{\|e\|_E\leq 1}} \left|\langle e,\phi \rangle \right| \lesssim \sup_{\genfrac{}{}{0pt}{}{e\in E}{\|e\|_{L^{t,s}}\leq 1}} \left|\langle e,\phi \rangle \right| \lesssim \|\phi\|_{F^*}.$$
These two inequalities with the properties of $\phi$ are impossible. Therefore we deduce that $F=L^{t,s}$ and so that $E$ is dense into $L^{t,s}$. Finally (\ref{propdual}) gives us that $E$ is equal to $L^{t,s}$ with equivalent norms (due to the completeness of the spaces).
\findem

\mb Make a mixture with Theorem \ref{complet}, then we finally obtain the main result~:

\begin{thm} \label{thm:inter2} Assume that $\mu(X)=\infty$ and the Hardy space $H^1$ satisfies ~: 
$$ H^1 \hookrightarrow L^1.$$
Assume that for $\sigma> \B'$, the operator $M_\sigma$ is bounded by $M_{HL,\B'}$.
Then for all $\theta\in (0,1),s\in[1,\infty)$, for all exponents $p\in(1,\B]$ and $t\in(1,\infty)$ such that 
$$ \frac{1}{t}= \frac{1-\theta}{p} + \theta <\frac{1}{\sigma'},$$
we have the equivalence between the two norms  
$$ \|\cdot \|_{L^{t,s}} \simeq \|\cdot \|_{(L^{p},H^1)_{\theta,s}}$$
and
$$ (L^{p},H^1)_{\theta,s}:= L^{t,s} \cap (H^1+L^p).$$
\end{thm}

\begin{rem} Note that the norm of the intermediate space does not depend on the Hardy space $H^1$ considered. 
\end{rem}

\begin{rem} We have seen in \cite{BJ} with the example of Riesz transforms, that the range of exponents where we can obtain Lebesgues spaces as intermediate spaces is optimal under the above assumption. 
\end{rem}

\subsection{Interpolation between Hardy spaces and weighted Lebesgue spaces.} \label{subsectionpoids}

We present in this subsection the weighted version of the previous results. For convenience, we assume that $\mu(X)=\infty$. We firstly recall the different class of weights~:

\begin{df}[The Muckenhoupt classes] A nonnegative function $\omega$ on $X$ belongs to the class ${\mathbb A}_p$ for $p\in(1,\infty)$ if
$$\sup_{Q \textrm{ ball}} \left(\frac{1}{\mu(Q)}\int_Q w d\mu \right) \left( \frac{1}{\mu(Q)}\int_Q \omega^{-1/(p-1)} d\mu \right)^{p-1} <\infty.$$
\end{df}

\begin{df}[The Reverse H\"older classes] A nonnegative function $\omega$ on $X$ belongs to the class $RH_q$ for $q\in(1,\infty)$, if there is a constant $C$ such that for every ball $Q\subset X$
$$ \left( \frac{1}{\mu(Q)}\int_Q \omega^q d\mu \right)^{1/q} \leq C \left(\frac{1}{\mu(Q)}\int_Q \omega d\mu \right).$$
\end{df}

\mb For the weight $\omega$, we define the associated measure (written by the same symbol) $\omega$ by $d\omega:=\omega d\mu$ and we denote $L^p_\omega:=L^p(X,d\omega)$ the corresponding weighted Lebesgue space. For the considered weights, the space $L^\infty_\omega$ does not depend on $\omega$ and will always be equal to $L^\infty$. 

\mb We have the well-known following properties (chapter 9 of \cite{Gra} for the Euclidean case)~:

\begin{prop} \label{prora} $1)$ For $s\in(1,\infty)$ the maximal operator $M_{HL,s}$ is bounded on $L^p(\omega)$ for all $p\in(s,\infty)$ and $\omega\in {\mathbb A}_{p/s}$. \\
$2)$ For $ p\in(1,\infty)$ and $\omega$ an ${\mathbb A}_p$-weight, there exists some constants $C,\epsilon>0$ such that for all balls $Q$ and all measurable subsets $A\subset Q$, we have~:
 \be{poidsh} \frac{\omega(A)}{\omega(Q)} \leq C \left(\frac{\mu(A)}{\mu(Q)} \right)^{\epsilon}. \ee
$3)$ For $\omega$ a nonnegative function and $p\in(1,\infty)$, we have the following equivalence~:
$$ \omega \in{\mathbb{A}}_{p} \Longleftrightarrow \omega^{1-p'} \in {\mathbb{A}}_{p'}.$$
\end{prop} 

\gb We give the weighted version of Theorem \ref{thm:interpolation}~:

\begin{thm} \label{thm:interpolationw} Assume that for $\sigma> \B'$, the operator $M_\sigma$ is bounded by $M_{HL,\B'}$. 
Let $\theta\in (0,1]$, $s\in[1,\infty]$, $q\in[\B',\infty)$ and $t$ satisfying
$$ \frac{1-\theta}{q} = \frac{1}{t}>\frac{1}{\sigma}.$$ 
Let $\omega$ be a weight belonging to $RH_{(\sigma/t)'}\cap {\mathbb A}_{q/\B'}$.
Then there exists a constant $c=c(\theta,s,\sigma,q)$ such that for all function $f\in L^{\B'}\cup L^q$~: if $f\in (L^q_\omega,(H^1)^*)_{\theta,s}$ then $f\in L^{t,s}_\omega$ and
$$ \|f\|_{L^{t,s}_\omega} \leq c \|f\|_{(L^q_\omega,(H^1)^*)_{\theta,s}}.$$
\end{thm}

\dem The proof is similar to the one of Theorem \ref{thm:interpolation}. It is based on Lemma \ref{good}. Using the same ideas and the weighted Theorem 3.1 of \cite{AM}, we prove this following one~:
\begin{lem} \label{goodw} Assume that for $\sigma> \B'$, $M_\sigma$ is bounded by $M_{HL,\B'}$. Let for $s\in(1,\frac{\sigma}{\B'}]$, $\omega$ be a weight belonging to $RH_{s'}$. Then for all $\rho$ large enough and for all $\gamma<1$ we have a ``good lambdas inequality''~: for all $\lambda>0$ 
$$ \omega\left( \left\{ M_{HL,\B'}(h) >\rho\lambda, M_{\B'}^\sharp(h)(x)\leq \gamma\lambda \right\} \right) \lesssim \left( \gamma \rho^{-\B'}+ \rho^{-\sigma} \right)^{1/s} \omega\left( \left\{ M_{HL,\B'}(h)>\lambda \right\} \right).$$
The implicit constant does not depend and $\rho$ and $\gamma$.
 In addition we have that for all exponent $q\in[\B',\sigma)$ there exists a constant $c=C(q)$ such that for all function $h\in L^{\B'}$
\be{equivw}  c\|h\|_{L^q_\omega} \leq \|M_{\B'}^\sharp(h) \|_{L^q_\omega} \leq c^{-1} \|h\|_{L^q_\omega}. \ee
\end{lem}

\mb Then we obtain the weighted version of Proposition \ref{thm:bonlambda}~:

\begin{prop} \label{inter2ww} Assume that for $\sigma> \B'$, the operator $M_\sigma$ is bounded by $M_{HL,\B'}$. Then for $q\in[\B',\infty)$, $s\in(1, \frac{\sigma}{\B'}]$ and a weight $\omega\in RH_{s'} \cap {\mathbb A}_{q/\B'}$, for all $\kappa\geq 2$ there is a constant $c_2=c_2(q,\kappa)\geq 1$ such that for all $t>0$ we have the following estimate of the $K$-functional~: for every function $f\in L^{\B'}$
\begin{align} 
\lefteqn{ K(t,M_{HL,q}(f),L^{q,\infty}_\omega,L^\infty) \leq c_2 \Big[K(t,f,L^q_\omega,(H^1)^*)  } & & \nonumber \\
 & &   + \kappa^{1-\sigma/(sq)}  K(\kappa^{\sigma/(sq)} t,M_{HL,q}(f),L^{q,\infty}_\omega,L^\infty) \Big]. \label{ineqK2w}
\end{align}
\end{prop}

\dem The proof is exactly the same as the one of Theorem \ref{thm:bonlambda}. We just have to check that the maximal operator $M^\sharp_{\B'}$ is of weak type $(q,q)$ for the new measure $\omega$. Under our assumption, this operator is bounded by $M_{HL,\B'}$, which is well bounded on $L^q_\omega$ as $\omega \in {\mathbb A}_{q/\B'}$. \findem 

\mb Then the end of the proof is the same, and so Theorem \ref{thm:interpolationw} is proved (by taking $s=\sigma/t$). We obtain also the following result~:

\begin{thm} \label{thm:interw} Assume that the Hardy space $H^1$ is complete, $\mu(X)=\infty$ and that for $\sigma> \B'$, the operator $M_\sigma$ is bounded by $M_{HL,\B'}$. Let $p\in(1,\B]$, $\theta\in (0,1),s\in[1,\infty)$ and $t\in(1,\infty)$ satisfying
$$ \frac{1}{t}= \frac{1-\theta}{p} + \theta <\frac{1}{\sigma'}.$$
Let $\omega$ be a weight in $RH_{\frac{1-t'}{1-\frac{p'}{b'}}} \cap {\mathbb A}_{t/\sigma'}$ and assume that $H^1$ is continuously embedded into $L^1_\omega$. Then we have the following equivalence~: 
\be{injectionw} L^{t,s}_\omega \simeq (L^{p}_\omega,H^1)_{\theta,s}.\ee
\end{thm}

\dem We let the detailled proof to the reader. The proof is analog to the one of Theorem \ref{thm:inter} using in addition some weighted arguments. We use the fact that (see Proposition 2.1 of \cite{AM}) : for $1\leq x\leq \infty$ and $1\leq y<\infty$
\be{poidseq}
\omega \in {\mathbb A}_{x} \cap RH_{y} \Longleftrightarrow \omega^{y}\in {\mathbb A}_{y(x-1)+1}.
\ee
So we can compute that
$$\omega \in RH_{\frac{1-t'}{1-\frac{p'}{b'}}} \cap {\mathbb A}_{t/\sigma'} \Longleftrightarrow \omega^{1-t'} \in RH_{(\sigma/t')'}\cap {\mathbb A}_{p'/\B'}.$$
Then we use Theorem \ref{thm:interpolationw} with the weight $\omega^{1-t'}$. We have to be careful because we are using duality with respect to the measure $\mu$ therefore 
$$ \left(L^{t,s}_\omega\right)^{*} = L^{t',s'}_{\omega^{1-t'}}.$$
\findem

\begin{rem} The constant appearing in (\ref{injectionw}) for the embedding
 \be{impor} L^{t,s}_\omega \hookrightarrow (L^{p}_\omega,H^1)_{\theta,s} \ee
is independent on the constant of the embedding $H^1 \hookrightarrow L^1_{\omega}$.
So to check this last inclusion, we can also assume that the weight is bounded and that $H^1 \hookrightarrow L^1$. Then the implicit constant in (\ref{impor}) does not depend on $\|\omega\|_\infty$.
\end{rem}

\mb We introduce the following definition due to \cite{AM}~:

\begin{df} For $\omega$ a nonnegative function on $X$ and $0<p_0<q_0\leq \infty$ two exponents, we introduce the set
$$ \mathcal{W}_{\omega}(p_0,q_0):=\left\{ p\in(p_0,q_0),\ \omega\in {\mathbb A}_{p/p_0} \cap RH_{(q_0/p)'} \right\}.$$
\end{df}

\begin{rem} \label{rempoids}  With the notation of Theorem \ref{thm:interw}, since $t'\leq p'$ then $\frac{1-t'}{1-\frac{p'}{\B'}} \geq (\frac{\B}{t})'$.
Thus if $\omega$ is a weight in $RH_{\frac{1-t'}{1-\frac{p'}{\B'}}} \cap {\mathbb A}_{t/\sigma'}$, then $t\in \mathcal{W}_{\omega}(\sigma',p)$.
\end{rem}

\section{Applications.}

\label{sectionapp}

\subsection{The linear theory.} \label{linear}

\gb 
To apply the previous abstract results, it is important to know when an operator is continuously acting in our Hardy spaces. We recall the following result of \cite{BJ}~:

\begin{prop} \label{prop:cont} Let $T$ a linear operator satisfying the following ``off-diagonal'' estimates~: for all ball $Q$, for all $j\geq 0$ there exist coefficients $\gamma_j$ such that for all $L^\B$-functions $f$ supported in $Q$
\be{hyp2h} \left(\frac{1}{\mu(2^{j+1}Q)} \int_{S_j(Q)} \left| T(B_Q(f))\right| d\mu \right) \leq \gamma_j\frac{\mu(Q)}{\mu(2^jQ)} \left(\frac{1}{\mu(Q)}\int_{Q} |f|^\B d\mu \right)^{1/\B}.\ee
If the coefficients $\gamma_j$ satisfy
\be{assum} \Lambda:= \sum_{j\geq 0} \gamma_j <\infty, \ee
then there exists a constant $C=C(\Lambda)$ such that for all atom $f\in H^1_{ato}$
$$ \|T(f)\|_{L^1} \leq C.$$
Consequently $T$ is continuous from $H^1_{F,ato}$ into $L^1$.
\end{prop}

\mb In \cite{BJ}, the authors have given a molecular version too.

\begin{rem} 
\mb $1-)$ The coefficients $\gamma_j$ may depend on the ball $\gamma_j=\gamma_j(Q)$ and then we have to replace (\ref{assum}) by
$$  \sup_{Q\ \textrm{ball}} \sum_{j\geq 0} \gamma_j(Q) <\infty.$$
$2-)$ We can just assume that $T$ is a sublinear operator or a positive linearizable operator which is meaning that there exists a Banach space $\mathcal{B}$ and a linear operator $U$ defined from $L^\B$ into $L^\B(X,\mathcal{B})$ such that
$$\forall f\in L^\B, \qquad T(f)(x) = \left\| U(f)(x)\right\|_{\mathcal{B}}, \qquad \mu-a.e.\ .$$
This improvement is useful to study some maximal operators.
\end{rem}

\mb We would like to deduce that $T$ can be continuously extended on the whole Hardy space $H^1_{ato}$. We know from the work \cite{Bownik} of M. Bownik that it is not sufficient in the general case to have boundedness on all the atoms. Using ideas of \cite{BJ}, we can found an operator $U$ bounded from $H^1_{ato}$ into $L^1$, which coincides with $T$ on $H^1_{F,ato}$. However we do not know if $U$ and $T$ coincide on more general functions in the Hardy space. \\
Fortunately to use interpolation, it is not a problem. So we will describe our interpolation result and we will finish this section by giving some conditions that permit us to extend our operator to the whole Hardy space.

\begin{thm} \label{theo4} Let us assume $1\leq \sigma'<p_0\leq \B$, $H^1_{ato}\hookrightarrow L^1$ and $\mu(X)=\infty$. Let $T$ be an $L^{p_0}$-bounded sublinear operator such that for all balls $Q$ and for all functions $f$ supported in $Q$
\be{hypp1} \forall j\geq 0 \quad \left(\frac{1}{\mu(2^{j+1}Q)} \int_{S_j(Q)} \left| T(B_Q(f))\right|^{p_0} d\mu \right)^{1/p_0} \leq \gamma_j\frac{\mu(Q)}{\mu(2^jQ)} \left(\frac{1}{\mu(Q)}\int_{Q} |f|^\B  d\mu \right)^{1/\B} \ee
and
$$\forall j\geq 0 \quad \left(\frac{1}{\mu(2^{j+1}Q)} \int_{S_j(Q)} \left| f-B_Q(f)\right|^\B d\mu \right)^{1/\B} \leq \gamma_j\frac{\mu(Q)}{\mu(2^jQ)} \left(\frac{1}{\mu(Q)}\int_{Q} |f|^{\sigma'} d\mu \right)^{1/\sigma'},$$
where the coefficients $\gamma_j$ satisfy
$$ \sum_{j\geq 0} \gamma_j <\infty.$$
Then for all exponents $p\in (\sigma',p_0)$, there exists a constant $C=C(p)$ such that
$$ \forall f\in L^{p_0} \cap L^p, \qquad  \|T(f)\|_{L^p} \leq C \|f\|_{L^p}.$$
\end{thm}

\begin{rem} We have seen in \cite{BJ} that for this particular application, we can choose a space $X$ of finite measure then the extra term is not a problem. \\
In \cite{BJ}, we have already proved this result without the assumption $H^1_{ato}\hookrightarrow L^1$ for linearizable operators. Here we want to explain the proof using real interpolation theory and so the result is true and new for general sublinear operators.
\end{rem}

\dem In \cite{BJ}, we shew that under these assumptions, the maximal
operator $M_{\sigma}$ (defined by (\ref{opeM})) is bounded by $M_{HL,\B'}$. In addition we know that $T$ is $L^p$-bounded. From Proposition \ref{prop:cont} we know that $T$ is bounded from $H^1_{F,ato}$ into $L^1$. We can use real interpolation theory and obtain for $\theta\in (0,1)$ and $s\in (1,\infty)$ the boundedness of $T$ from $(L^p,H^1_{F,ato})_{\theta,s}$ to $(L^p,L^1)_{\theta,s}$. We know that this last intermediate space corresponds to the Lebesgue space. In addition from Theorem \ref{thm:inter2}, we know that the norm $\|\ \|_{(L^p,\overline{H^1_{F,ato}})_{\theta,s}}$ is equivalent to the norm in $L^{t,s}$ and that
$$ (L^p,H^1_{F,ato})_{\theta,s} = L^{t,s} \cap L^p,$$
with 
$$\frac{1}{t}=\theta+ \frac{1-\theta}{p_0}<\frac{1}{\sigma'}.$$
We then obtain the desired conclusion with $t=p=s$. \findem

\begin{rem} \label{remi} In \cite{BJ}, we have already obtained a weighted result for linearizable operators satifying (\ref{hypp1}) : consider $T$ a linearizable operator of Theorem \ref{theo4} then for all weight $\omega$ and every exponent $p\in {\mathcal W}_{\omega}(\sigma',\B)$, there is a constant $C$ such that
$$ \forall f\in L^{p_0} \cap L^p_\omega, \qquad  \|T(f)\|_{L^p_\omega} \leq C \|f\|_{L^p_\omega}.$$
With our previous weighted results, we can describe a similar weighted version for a general sublinear operator but we have to require 
$$\omega \in RH_{\frac{1-p'}{1-\frac{p_0'}{\B'}}} \cap {\mathbb A}_{p/\sigma'},$$
which is a stronger assumption due to Remark \ref{rempoids}. However it is interesting to note that this stronger assumption permit us to obtain a weighted interpolation result for sublinear operators, where as the weaker condition ($p\in {\mathcal W}_{\omega}(\sigma',\B)$) requires a linearizable operator due to the use of duality.
\end{rem}

\mb We finish this subsection by studying the following problem~: let $T$ be an operator bounded on all the atoms in the space $L^1$. Can we extend it continuously in the whole Hardy space $H^1_{ato}$ with its natural norm ? This problem was studied for particular case in several papers (see \cite{MSV} and \cite{YZ}).
Following the ideas of \cite{MSV}, we have the following result~:

\begin{prop} Assume that the Hardy space $H^1_{ato} \hookrightarrow L^1$ and that 
$$ \bigcap_{Q\textrm{ balls}} \textrm{ker} (B_Q^*) \subset L^\infty.$$
Here we use the notation of Lemma \ref{lem1} : $B_Q^*$ is acting on $L^{\B'}+(H^1_{ato})^*$. \\
Let $T$ be a $L^\B$-bounded linear operator with a constant $C$ such that for all atoms $f\in H^1_{ato}$, we have
 $$ \|T(f)\|_{L^1} \leq C.$$
Then it can be continuously extended on $H^1_{ato}$ into $L^1$.
\end{prop}

\dem We know (see \cite{BJ} or \cite{MSV}) that there exists an operator $U$ continuous from $H^1_{ato}$ into $L^1$ such that for each atom $m$ : $U(m)=T(m)$. We have to prove that
$$ \forall f\in L^\B \cap H^1_{ato}, \qquad U(f)=T(f).$$
To prove this fact, we use duality. Let $f\in L^\B$ compactly supported and $g\in L^\infty \cap L^{\B'}$. For all balls $Q$, $B_Q(f)\in H^1_{F,ato}$ because $\mu$ is a borelian measure. So
$$ \langle T(B_Qf),g\rangle = \langle U(B_Q(f)), g\rangle.$$
We deduce that
$$ \langle f,B_Q^{*}T^*g\rangle = \langle f,B_Q^*U^* g\rangle.$$
Hence for all compactly supported functions $f\in L^\B$,
$$\langle f, B_Q^*(T^*g- U^*g) \rangle=0,$$
therefore $B_Q^*(T^*g- U^*g)=0$. We know that $h:=T^*g- U^*g\in L^{\B'}+(H^1_{ato})^*$.
Under our assumption, we conclude that $h\in L^\infty$ and so for each atom $m$
$$ \langle m, h \rangle =0.$$
 As the Hardy space $H^1_{ato}$ is embedded into $L^1$, by the duality $L^1-L^\infty$, we know that 
$\langle m, h \rangle=\int mh d\mu$. Thus for all functions $f\in H^1_{ato}$ we have 
$$ \langle f, h \rangle=0.$$
In particular for $f\in L^\B \cap H^1_{ato}$, we get
$$ \langle f, h \rangle=0 = _{L^{\B}}\langle f,T^*g\rangle_{L^{\B'}} - _{H^1_{ato}}\langle f,U^*g\rangle_{(H^{1}_{ato})^*} = _{L^{\B}}\langle T(f),g\rangle_{L^{\B'}} - _{H^1_{ato}}\langle U(f),g\rangle_{(H^{1}_{ato})^*}.$$
This is true for all functions $g\in L^\infty \cap L^{\B'}$. We deduce that $T(f)=U(f)$ in $\left(L^\infty \cap L^{\B'}\right)^*$ and so $T(f)(x)=U(f)(x)$ for almost every $x\in X$. \findem

\subsection{The bilinear theory.}
\label{bilinear}

In all this subsection, we implicitly use a space of homogeneous type $X$ with an infinite measure : $\mu(X)=\infty$, in order to use our interpolation result of Theorem \ref{thm:inter2}.

\gb We are interested in a bilinear version of Theorem \ref{theo4}. We choose two collections ${\mathbb B}^1:=(B_Q^1)_{Q\in {\mathcal Q}}$ and ${\mathbb B}^2:=(B_Q^2)_{Q\in {\mathcal Q}}$ and two exponents $\B_1,\B_2\in(1,\infty]$. We assume that 
${\mathbb B}^1$ is a collection of $L^{\B_1}$-bounded operators and that ${\mathbb B}^2$ is a collection of $L^{\B_2}$-bounded operators. We can also define two kinds of Hardy spaces $H^1_{{\mathbb B}^1}$ and $H^1_{{\mathbb B}^2}$. According to Definition \ref{hardyF}, we can construct the spaces $H^1_{F,ato,{\mathbb B}^1}$ and $H^1_{F,ato,{\mathbb B}^2}$.
In this context, we have also the following bilinear results~:

\begin{prop} \label{prop:bilineaire} Let $T$ be a bilinear operator with coefficients $(\gamma_j)_{j\geq 0}$ satisfying for all balls $Q_1,Q_2$ and for all functions $f,g$ supported in $Q_1$ and $Q_2$, we have for $l=1,2$ and for all $j_1,j_2\geq 0$ 
\begin{align} 
 \lefteqn{\left(\frac{1}{\mu(2^{j_l+1}Q_l)} \int_{C_{j_1}(Q_1)\cap C_{j_2}(Q_2)} \left| T(B^1_{Q_1}(f),B^2_{Q_2}(g))\right| d\mu \right) \lesssim} & & \nonumber \\ 
 & & \gamma_{j_1} \gamma_{j_2} \frac{\mu(Q_1)}{\mu(2^{j_1}Q_1)} \frac{\mu(Q_2)}{\mu(2^{j_2}Q_2)} \left(\frac{1}{\mu(Q_1)}\int_{Q_1} |f|^{\B_1}  d\mu \right)^{1/\B_1}\left(\frac{1}{\mu(Q_2)}\int_{Q_2} |g|^{\B_2}  d\mu \right)^{1/\B_2} \label{bi1} 
\end{align}
with coefficients $\gamma_l$ satisfying
\be{propcoeff} \sum_{l\geq 0} \gamma_l \lesssim 1. \ee
Then the operator $T$ is continuous from $H^1_{F,ato,{\mathbb B}^1} \times H^1_{F,ato,{\mathbb B}^2}$ into $L^{1/2}$.
\end{prop}

\dem We will use ideas of \cite{GK}. Let $f\in H^1_{F,ato,{\mathbb B}^1}$ and $g\in H^1_{F,ato,{\mathbb B}^2}$, we can also write them with a finite atomic decomposition :
$$f=\sum_{Q} \lambda_Q B_Q^1(f_Q) \qquad g=\sum_{R} \tau_R B_R^2(g_R)$$
with the appropriate properties for $f_Q$ and $g_R$~: $f_Q$ is supported in $Q$ with $\|f_Q\|_{\B} \leq \mu(Q)^{-1/\B'}$, 
\be{propcoeff2} \sum_{Q} |\lambda_Q| \leq 2 \|f\|_{H^1_{{\mathbb B}^1}} \ee
and similarly for $g_R$, relatively to the ball $R$.
So we have to study
$$T(f,g)=\sum_{Q,R} \lambda_Q \tau_R T(B_Q^1f_Q,B_R^2g_R).$$
We decompose with the coronas around the balls 
$$T(f,g)=\sum_{Q,R} \sum_{i,j\geq 0} \lambda_Q \tau_R T(B_Q^1f_Q,B_R^2g_R){\bf 1}_{C_i(Q)\cap C_j(R)}.$$
To estimate the norm $\|T(f,g)\|_{L^{1/2}}$ by symmetry we need just to study the sum over the extra condition
\be{extracond} 2^ir_Q \leq 2^jr_R,\ee
where $r$ corresponds to the radius of the ball. We recall Lemma 2.1 of \cite{GK}~:
\begin{lem} For $r\leq 1$, there exists a constant $C$ and $\delta>1$ such that for all collection $(Q_k)_k$ of balls and $(g_k)_k$ collection of nonnegative integrable functions supported in $Q_k$ we have
$$\left\| \sum_{k} g_k \right\|_{L^r} \leq C \left\| \sum_{k} \left(\frac{1}{\mu(Q_k)}\int g_k\right){\bf 1}_{\delta Q_k} \right\|_{L^r}.$$ 
\end{lem}

\dem The proof is explained in the particular case of $\R^n$, using the dyadic structure of the euclidean space. However the proof can easily be extended to a general space of homogeneous type using its dyadic structure (proved in \cite{David}) so we let the details to the reader. \findem

\mb Now using this Lemma with $r=1/2$, we get
\begin{align*}
\left\|T(f,g)\right\|_{L^{1/2}} & \lesssim \left\| \sum_{Q,R} \sum_{i,j\geq 0} \left|\lambda_Q\right| \left|\tau_R\right| \left|T(B_Q^1f_Q,B_R^2g_R)\right| {\bf 1}_{C_i(Q)\cap C_j(R)} \right\|_{L^{1/2}} \\
 & \lesssim \left\| \sum_{Q,i} \left|\lambda_Q\right| \left(\sum_{\genfrac{}{}{0pt}{}{R,j}{}} \left|\tau_R\right| \left|T(B_Q^1f_Q,B_R^2g_R)\right|{\bf 1}_{C_i(Q)\cap C_j(R)}\right) {\bf 1}_{2^{i}Q} \right\|_{L^{1/2}} \\
 & \lesssim \left\| \sum_{Q,R} \sum_{i,j\geq 0} \left|\lambda_Q\right| \left|\tau_R\right|  \left( \frac{1}{\mu(2^iQ)}\int_{C_i(Q)\cap C_j(R)} \left|T(B_Q^1f_Q,B_R^2g_R)\right| d\mu \right) {\bf 1}_{\delta 2^{i}Q} \right\|_{L^{1/2}}.
\end{align*}
We can use the doubling property of the measure, the estimates (\ref{bi1}) and (\ref{extracond}) to finally obtain
\begin{align*}
 \lefteqn{\left\|T(f,g)\right\|_{L^{1/2}} \lesssim \left\| \sum_{Q,R} \sum_{i,j\geq 0} \left|\lambda_Q\right| \left|\tau_R\right| \gamma_{i} \gamma_{j} \right. } & & \\
 & & \left. \frac{\mu(Q) \mu(R)}{\mu(2^{i}Q)\mu(2^{j}R)} \left(\frac{1}{\mu(Q)}\int_{Q} |f_Q|^{\B_1}  d\mu \right)^{1/\B_1}\left(\frac{1}{\mu(R)}\int_{R} |g_R|^{\B_2}  d\mu \right)^{1/\B_2} {\bf 1}_{\delta 2^{i}Q} {\bf 1}_{2\delta 2^{j}R} \right\|_{L^{1/2}}. 
\end{align*}
We can ``add'' the function ${\bf 1}_{2\delta 2^{j}R}$ due to the condition (\ref{extracond}).
With the estimates on $f_Q$ and $g_R$, we have
 \begin{align*}
 \left\|T(f,g)\right\|_{L^{1/2}} &  \lesssim \left\| \sum_{Q,R} \sum_{i,j\geq 0} \left|\lambda_Q\right| \left|\tau_R\right|  
 \frac{\gamma_{i} \gamma_{j}}{\mu(2^{i}Q)\mu(2^{j}R)} {\bf 1}_{\delta 2^{i}Q} {\bf 1}_{2\delta 2^{j}R} \right\|_{L^{1/2}} \\
& \lesssim \left\| \sum_{Q,R} \sum_{i,j\geq 0} \left|\lambda_Q\right| \left|\tau_R\right|  
 \frac{\gamma_{i}}{\mu(2^iQ)}{\bf 1}_{\delta 2^{i}Q} \frac{\gamma_{j}}{\mu(2^{j}R)} {\bf 1}_{2\delta 2^{j}R} \right\|_{L^{1/2}}.
\end{align*}
By using H\"older inequality, we get
 \begin{align*}
 \left\|T(f,g)\right\|_{L^{1/2}} \lesssim \left\| \sum_{Q,i} \left|\lambda_Q\right|
 \frac{\gamma_{i}}{\mu(2^iQ)}{\bf 1}_{\delta 2^{i}Q} \right\|_{L^1} \left\| \sum_{R,j} \left|\tau_R\right|\frac{\gamma_{j}}{\mu(2^{j}R)} {\bf 1}_{2\delta 2^{j}R} \right\|_{L^1}.
\end{align*}
Then the proof is finished with the properties (\ref{propcoeff}) and (\ref{propcoeff2}) and the doubling property of the measure $\mu$.
\findem

\mb If we would like to use the Hardy space for just one of the two functions, we have the following version~:

\begin{prop} \label{prop:bilineaire2} Let assume $1<\B_1,\B_2$. Let $T$ be a bilinear operator with coefficients $(\gamma_j)_{j\geq 0}$ such that for all balls $Q_1$ and all balls $Q_2$ with radius $r_{Q_2}=1$, for all functions $f,g$ supported in $Q_1$ and $Q_2$, we have for $l=1,2$ and for all $j_1,j_2\geq 0$ 
\begin{align}
\lefteqn{ \left(\frac{1}{\mu(2^{j_l+1}Q_l)} \int_{C_{j_1}(Q_1)\cap C_{j_2}(R_2)} \left| T(B_{Q_1}f,g)\right| d\mu \right) \lesssim} \nonumber & & \\
 & & \gamma_{j_1} \gamma_{j_2} \frac{\mu(Q_1)}{\mu(2^{j_1}Q_1)} \frac{\mu(Q_2)}{\mu(2^{j_2}Q_2)} \left(\frac{1}{\mu(Q_1)}\int_{Q_1} |f|^{\B_1}  d\mu \right)^{1/\B_1} \left(\frac{1}{\mu(Q_2)}\int_{Q_2} |g|^{\B_2}  d\mu \right)^{1/\B_2} \label{bi10} 
\end{align}
with coefficients $\gamma_l$ satisfying
\be{propcoeff20} \sum_{l\geq 0} \gamma_l \lesssim 1. \ee
Then the operator $T$ is continuous from $H^1_{F,ato,{\mathbb B}^1} \times L^{\T_2}$ into $L^{\T}$ for all exponents satisfying $0<\B_2<\T_2 \leq \infty$ and
\be{ex} \frac{1}{\T}=\frac{1}{\T_2}+1.\ee
\end{prop}

\begin{rem} 
The scale $r_{Q_2}=1$, taken for the balls $Q_2$ is not important. The important fact is to have these ``off-diagonal'' decays for all balls $Q_2$ at a fixed scale.
 \end{rem}

\begin{rem} \label{remim} Let assume that $T$ satisfies Proposition \ref{prop:bilineaire2} with $\B_2=\infty$ then, by the same arguments as we will explain in the proof, we can show that for each function $g\in L^\infty$, the linear operator $f\rightarrow T(f,g)$ satisfies Theorem \ref{theo4} about linear operators.
\end{rem}

\dem The proof is similar to the previous one with choosing $B^2_{Q_2}=Id$. 
First we will prove a weak type estimate for the bilinear operator. So let $R$ be a set of finite measure. We fix a function $g$ supported on $R$ and bounded by $1$. We use a bounded covering of the space $X$ with balls of radius $1$ : $(Q_2)_{Q_2}$. We decompose the function $g$ over this covering and so as previously with $\T_2\geq \B_2$ we get~:
\begin{align*}
 \left\|T(f,g)\right\|_{L^\T} & \lesssim \left\| \sum_{Q_1,Q_2} \sum_{i,j\geq 0}   
 \frac{\gamma_{i}}{\mu(2^iQ_1)} {\bf 1}_{\delta 2^{i}Q_1} \frac{\gamma_{j}\mu(Q_2)}{\mu(2^{j}Q_2)} \left(\frac{\mu(R\cap Q_2)}{\mu(Q_2)}\right)^{1/\T_2} {\bf 1}_{2\delta 2^{j}Q_2} \right\|_{L^\T}.
\end{align*}
By H\"older inequality and property (\ref{propcoeff2}), as $\T_2\geq 1$ we obtain
\begin{align*}
 \left\|T(f,g)\right\|_{L^\T} \lesssim \|f\|_{H^1_{F,ato}} \sum_{j\geq 0} \gamma_j \left\| \sum_{Q_2} \frac{\mu(Q_2)}{\mu(2^{j}Q_2)} \left(\frac{\mu(R\cap Q_2)}{\mu(Q_2)}\right)^{1/\T_2} {\bf 1}_{2\delta 2^{j}Q_2} \right\|_{L^{\T_2}}.
\end{align*}
However there exists a constant $C=C(X)$ such that for all $x\in X$, for all $j\geq 0$
\be{astuce} \sum_{\genfrac{}{}{0pt}{}{Q_2}{x\in 2^{j}Q_2}} \frac{\mu(Q_2)}{\mu(2^{j}Q_2)} \leq C. \ee
This is due to the following fact : all the balls $Q_2$ have the same radius so for $Q_2^1,Q_2^2$ two balls considered in the sum, we have $\mu(2^jQ_2^1) \simeq \mu(2^jQ_2^2)$. Then we use that the collection of the balls $Q_2$ is a bounded covering.
Using (\ref{astuce}), we obtain
\begin{align*} 
\left\| \sum_{Q_2} \frac{\mu(Q_2)}{\mu(2^{j}Q_2)} \left(\frac{\mu(R\cap Q_2)}{\mu(Q_2)}\right)^{1/\T_2} {\bf 1}_{2\delta 2^{j}Q_2} \right\|_{L^{\T_2}}^{\T_2} & \lesssim 
 \int \sum_{Q_2}  \left(\frac{\mu(R_2\cap Q_2)}{\mu(Q_2)}\right) \frac{\mu(Q_2)}{\mu(2^{j}Q_2)}{\bf 1}_{2\delta 2^{j}Q_2} d\mu  \\
& \lesssim \sum_{Q_2} \mu(R\cap Q_2) \lesssim \mu(R).
\end{align*}
We have (again) used at the last inequality the fact that the collection of the balls is a bounded covering.
We deduce also that
$$\left\|T(f,g)\right\|_{\T} \lesssim \|f\|_{H^1_{F,ato}} \mu(R_2)^{1/\T_2}.$$
Then for $f$ fixed, the linear operator $T(f,.)$ is of weak type $(\T_2,\T)$ for all exponents $(\T_2,\T)$ satisfying (\ref{ex}). The strong continuities are obtained by real interpolation for Lorentz spaces.
\findem

\mb These two propositions will permit us to refind results about bilinear Calder\'on-Zygmund and related operators. We will explain these examples in the following section. 

\mb For convenience, we will use the following definition~:
\begin{df} Let $p_0\leq \B$ be exponents and $H^1$ be a Hardy space defined by a collection ${\mathbb B}$ of $L^\B$-bounded linear operators. We say that the space $H^1$ is ``$L^{p_0}-L^\B$ regularizing'' if the maximal operator $M_{p_0'}$ (defined by \ref{opeM}) is bounded by the Hardy-Littlewood operator $M_{HL,\B'}$ and if the Hardy space is embedded into $L^1$.
\end{df}

\begin{rem} \label{remms} Under the previous notations, we have seen that if $H^1 \hookrightarrow L^1$ there exist coefficients $\gamma_j$ such that for all balls $Q$ and all functions $f$ supported on $Q$
$$\forall j\geq 0 \quad \left(\frac{1}{\mu(2^{j+1}Q)} \int_{S_j(Q)} \left| f-B_Q(f)\right|^\B d\mu \right)^{1/\B} \leq \gamma_j\frac{\mu(Q)}{\mu(2^jQ)} \left(\frac{1}{\mu(Q)}\int_{Q} |f|^{p_0} d\mu \right)^{1/p_0},$$
satisfying
$$ \sum_{j\geq 0} \gamma_j <\infty,$$
then the Hardy space is $L^{p_0}-L^\B$ regularizing.
\end{rem}

\mb Now using bilinear interpolation, we obtain the following theorems.

\begin{thm} \label{thm:bilineaire1} Let assume that the Hardy space $H^1_{F,ato,{\mathbb B}_i}$ is  
$L^{q_i}-L^{\B_i}$ regularizing for $i=1,2$ and for exponents $1\leq q_i<\B_i$.
Let $T$ be a bilinear operator satisfying the assumptions of Proposition \ref{prop:bilineaire} and bounded from $L^{p_1} \times L^{p_2}$ into $L^{p}$ for exponents satisfying $p_i\in (q_i,\B_i]$ and
$$0<\frac{1}{p}=\frac{1}{p_1}+\frac{1}{p_2}\leq 1.$$
Then for all $\theta\in (0,1)$ such that
\begin{align*}
\frac{1}{\T_1} & :=\frac{1-\theta}{p_1}+\theta<\frac{1}{q_1} \\
\frac{1}{\T_2} & :=\frac{1-\theta}{p_2}+\theta<\frac{1}{q_2}
\end{align*}
there exists a constant $C=C(\theta)$ satisfying
$$ \forall f_i\in L^{\T_i} \cap L^{p_i}, \qquad  \|T(f_1,f_2)\|_{L^\T} \leq C \|f_1\|_{L^{\T_1}}\|f_2\|_{L^{\T_2}},$$
with an exponent $\T$ given by
$$\frac{1}{\T_1}+\frac{1}{\T_2}:=\frac{1}{\T}.$$
\end{thm}

\mb We can summarize this result with Figure 1. 

\begin{figure}[!h]
 \centering
 \input{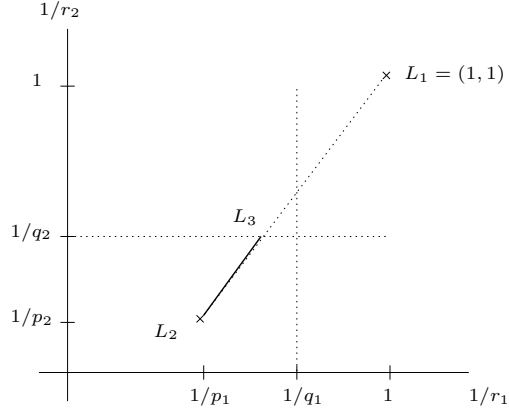}
 \label{figure1}
 \caption{Bilinear interpolation with the two Hardy spaces.}
\end{figure}

\mb So we have a strong continuity in the point $L_2$. Proposition \ref{prop:bilineaire} describes the continuity at the point $L_1$ and then by bilinear interpolation, we obtain strong continuities for all the points between $L_2$ and $L_3$.

\gb
\dem We have shown (in Proposition \ref{prop:bilineaire}) that the operator $T$ is continuous from $H^1_{{\mathbb B}^1} \times H^1_{{\mathbb B}^2}$ into $L^{1/2}$. We know by assumption that it is continuous from $L^{p_1} \times L^{p_2}$ into $L^{p}$ with $p,p_1,p_2\geq 1$. We use then bilinear real interpolation. We work with a Lorentz space $L^{1/2}$ (which  is a quasi-Banach space), the bilinear interpolation for these spaces is studied in \cite{GM} by L. Grafakos and M. Mastylo. By the Corollary 5.1 of their paper, we obtain that our operator $T$ is continuous from $(L^{p_1},H^1_{{\mathbb B}^1} )_{\theta,s_1} \times (L^{p_2},H^1_{{\mathbb B}^2})_{\theta,s_2}$ into the space $(L^{p},L^{1/2})_{\theta,s}$ for all $\theta\in(0,1)$ and all parameters $s_i,s\in[1,\infty]$ satisfying
$$\frac{1}{s}+1=\frac{1}{s_1}+\frac{1}{s_2}.$$
With our assumptions and taking $s_i=1$, Theorem \ref{thm:inter} gives us that 
$$(L^{p_i},H^1_{{\mathbb B}^i})_{\theta,s_i} \simeq L^{\T_i,s_i}=L^{\T_i,1}.$$
We already know that
$$(L^{p},L^{1/2})_{\theta,s} \simeq L^{\T,1}\hookrightarrow L^{\T,\infty}$$
So by bilinear interpolation, we know that $T$ is continuous from $L^{\T_1,1} \times L^{\T_2,1}$ into $L^{\T,\infty}$. We now use bilinear interpolation for Lorentz spaces (see \cite{GTa} and \cite{mtt2}) to conclude the proof and to obtain strong continuity from $L^{\T_1} \times L^{\T_2}$ into $L^{\T}$. \findem

\mb We can use interpolation on just one side and not on the two sides together.

\begin{thm} \label{thm:bilineaire2} Let assume that the Hardy space $H^1_{F,ato,{\mathbb B}^1}$ is  
$L^{q_1}-L^{\B_1}$ regularizing for an exponent $q_1<\B_1$.
Let $T$ be a bilinear operator satisfying the assumptions of Proposition \ref{prop:bilineaire2} and bounded from $L^{p_1} \times L^{p_2}$ into $L^p$ with
$$0<\frac{1}{p}=\frac{1}{p_1}+\frac{1}{p_2}\leq 1,$$
$p_2=\B_2$ and $p_1\in(q_1,\B_1)$.
Then for all exponent $(\T_1,\T_2,\T)$ such that
$$\frac{1}{\T}=\frac{1}{\T_1}+\frac{1}{\T_2}$$
and 
$$ 0<\frac{\frac{1}{p_2}-\frac{1}{\T_2}}{\frac{1}{p_1}-\frac{1}{\T_1}}<\frac{ \frac{1}{p_2}}{\frac{1}{p_1}-1}, $$
the operator $T$ admits a continuous extension from $L^{\T_1}\times L^{\T_2}$ into $L^{\T}$.
\end{thm}

\mb We can ``think'' this result with Figure 2.

\begin{figure}[!h]
 \centering
 \input{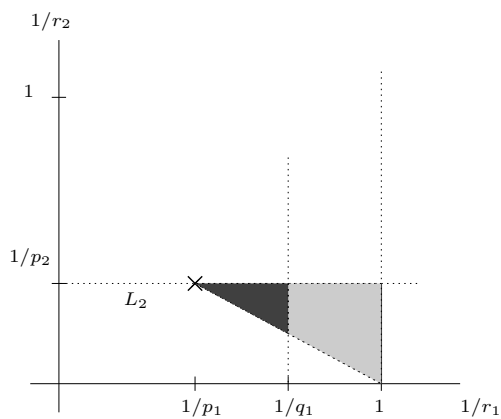}
 \label{figure2}
  \caption{Bilinear interpolation with one of the two Hardy spaces.}
\end{figure}

\dem Proposition \ref{prop:bilineaire2} shows us that we have a strong continuity for $T$ on the line $\T_1=1$ (with $\T_2\geq \B_2$). We assume that we have a first strong continuity on a point $(p_1,p_2)$ (which is the point $L_2$ on the figure with $p_2=\B_2$). By bilinear interpolation, we can deduce strong continuities in all the hatched region. We conclude the proof by using Theorem \ref{thm:inter} : we know that for $\T_1\geq q_1$ the intermediate space corresponds to the Lebesgue space. So in the most hatched region, we have strong continuities in Lebesgue space. This domain is exactly decribes by the condition over the exponents. \findem

\gb Then we can compute the two previous Theorems and combine them with duality and bilinear interpolation to get other results, we let that to the reader. Similarly, we could use interpolation with weighted Lebesgue spaces. For the moment, it is not clear what kind of results in weighted spaces should be reasonnable. That is why we decide to not describe the weighted version. However it is interesting to note that we could obtain continuities in $L^{\T_1}_{\omega_1} \times L^{\T_2}_{\omega_2}$ in $L^{\T}_\omega$ with different weights $\omega_1,\omega_2,\omega$. In addition, as we will see in the examples, it seems interesting to mixt these weighted results with the weighted extrapolation theory, described by L. Grafakos and J.M. Martell in \cite{GM2}.

\gb We would like to finish this section by an application of duality to obtain other bilinear continuities. We follow ideas of \cite{GT} and \cite{DGY}. In these two articles, the authors use duality and the end point estimate (the continuity from $H^1_{F,ato} \times H^1_{F,ato}$ into $L^{1/2}$) to obtain other continuities.

\begin{df} Let $T$ be a bilinear operator bounded from $L^{\B_1} \times L^{\B_2}$ into $L^\B$ for exponents $\B,\B_1,\B_2\in[1,\infty)$, we define its two adjoints $T^{*1}$ and $T^{*2}$ by
$$ \langle T(f,g) , h \rangle := \langle T^{*1}(h,g),f\rangle := \langle T^{*2}(f,h),g\rangle.$$
So $T^{*1}$ is bounded from $L^{\B'} \times L^{\B_2}$ into $L^{\B_1'}$ and $T^{*2}$ is bounded from $L^{\B_1} \times L^{\B'}$ into $L^{\B_2'}$. 
\end{df}

\begin{thm} \label{thm:final} Let assume that we have parameters $1\leq q<\B \leq \infty$ and a class $BS$ of bilinear operators such that~: 
\begin{enumerate}
 \item For $T\in BS$, its two adjoints $T^{*1},T^{*2}$ belong to the class $BS$. 
 \item For each operator $T\in BS$, there exist two Hardy spaces ($H^1_{F,ato,{\mathbb B}^1}$ and $H^1_{F,ato,{\mathbb B}^2}$) $L^q-L^\B$ regularizing such that $T$ is continuous from $H^1_{F,ato,{\mathbb B}^1} \times H^1_{F,ato,{\mathbb B}^2}$ into $L^{1/2}$. 
\end{enumerate}
Let $T$ be an operator of this class $BS$. Assume that there exists three exponents $p_1,p_2,p'\in(q,\B)$ satisfying
$$\frac{1}{p}=\frac{1}{p_1}+\frac{1}{p_2},$$
such that $T$ is continuous from $L^{p_1} \times L^{p_2}$ into $L^p$. Then for any exponents
$\T_1,\T_2,\T'\in(q,\B)$ satisfying
$$\frac{1}{\T}=\frac{1}{\T_1}+\frac{1}{\T_2}$$
the operator $T$ is continuous from $L^{\T_1} \times L^{\T_2}$ into $L^{\T}$.
\end{thm}

\dem The arguments are written in details in \cite{GT}. We describe them by the figures $3$ and $4$. In the first one (Figure 3), we have drawn the domain for the exponents $\T_1,\T_2$.

\begin{figure}[!h]
 \centering
 \input{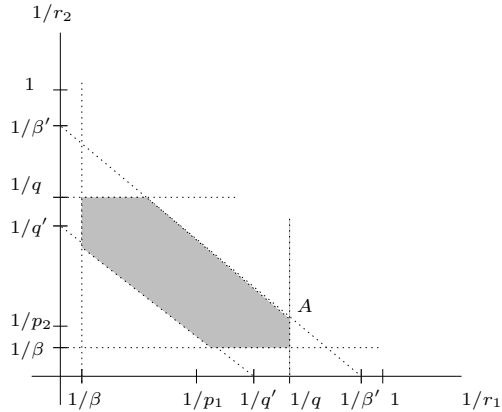}
 \caption{Domain of admissible exponents.}
 \label{figure3}
\end{figure}

\mb Now we will prove how to obtain the continuity near the point $A$, with the next figure.

\begin{figure}[!h]
 \centering
 \input{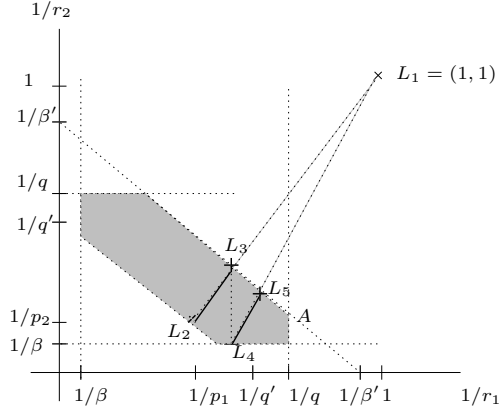}
 \caption{Scheme of the proof.}
 \label{figure4}
\end{figure}

\mb We start from the continuity on the point $L_2$. Then we use bilinear interpolation with the end point estimate at $L_1$ (this is described in Theorem \ref{thm:bilineaire1}). We obtain also continuity in Lebesgue space at the point $L_3$. Then we use duality to obtain continuity for $T^{*1}$ at the point $L_4$. By assumption $1-)$ of the class $BS$, we can repeat this procedure to $T^{*1}$ so we have continuity for it on the point $L_5$. By duality, we obtain also continuity for $T$ at a point $L_6$ (that we have not drawn for convenience) which is between the first point $L_2$ and $A$. By iterating this procedure, we can approach the point $A$, as closed as we want. By duality and symmetry, we obtain continuities for points as near as we want of the extremal points. Then by bilinear interpolation, we get continuities in whole the convex envelop, which corresponds to the hatched domain.
\findem

\gb All these results can seem to be technical. We will give examples in some well known cases to show how use them.

\section{Examples} \label{sectionexa}

\mb We remember that in \cite{BJ2}, we have already described an application of the linear theory for the problem of maximal $L^p$-regularity associated to a Cauchy problem. \\
In this section, we give examples and we explain how to use the previous bilinear abstract results.

\subsection{The bilinear Calder\'on-Zygmund operators.}

\mb Let $X$ be a space of homogeneous type with $\mu(X)=\infty$. Choose for our $B_Q$ operator, the exact oscillation
 $$B_Q(f):= f(x) - \left(\frac{1}{\mu(Q)}\int_Q f d\mu \right) {\bf 1}_{Q}.$$
In this case, we know that the Hardy space $H^1_{ato}$ corresponds to the classical Hardy space of Coifman-Weiss $H^1_{CW}$ (due to the work of E.M. Stein \cite{Se} and \cite{steine}).
Then the operator $A_Q$ is given by
$$ A_Q(f):=\left(\frac{1}{\mu(Q)}\int_Q f d\mu \right) {\bf 1}_{Q}$$ and we have that $A_Q=A_Q^*$.
We deduce also that
$$ M_{\infty}(f) \leq M_{HL,1}(f),$$
where $M_{\infty}$ is defined by (\ref{opeM}). So from Theorem \ref{thm:inter}, the Hardy space $H^1_{F,ato}$ is $L^{1}-L^{\infty}$ regularizing. 

\gb In this case, let $T$ be a bilinear operator associated to a bilinear kernel $k$ such that for all $f,g$ compactly supported and for all $x\in \textrm{supp}(f)^c \cup \textrm{supp}(g)^c$, we have the integral representation~:
$$T(f,g)(x):=\int_{X} \int_X k(x,y,z) f(y) g(z) d\mu(y) d\mu(z).$$
Let us assume that $k$ satisfy the ``standard'' bilinear estimates (for a certain $\epsilon>0$)~:
\be{cz1}
\left|k(x,y,z)-k(x,h,z)\right| \lesssim \left( \frac{d(y,h)}{ d(x,y)+d(z,x)} \right)^\epsilon \frac{1}{\mu(B(x,d(y,x))) + \mu(B(x,d(z,x)))^2} \ee for $d(y,h)\leq d(x,y)/2$ and
\be{cz2}
\left|k(x,y,z)-k(x,y,h)\right| \lesssim \left( \frac{d(z,h)}{ d(x,y)+d(z,x)} \right)^\epsilon\frac{1}{\mu(B(x,d(y,x))) + \mu(B(x,d(z,x)))^{2}} \ee
for $d(z,h)\leq d(x,z)/2$. 

\subsubsection{Continuities in Lebesgue spaces.}

\mb Suppose in addition that we already know a strong continuity for $T$ : $T$ is continuous from $L^{p_1}\times L^{p_2}$ to $L^{p}$ with $1\leq p_1,p_2,p<\infty$ and 
$$\frac{1}{p}=\frac{1}{p_1}+\frac{1}{p_2}.$$
Then we claim that this bilinear Calder\'on-Zygmund operator satisfies the required properties of Proposition \ref{prop:bilineaire} with $\B_1=\B_2=\infty$ and $H^1_{{\mathbb B}^1}=H^1_{{\mathbb B}^1}=H^1_{CW}$. We quickly check this fact. With the notations of this proposition : for $j_1=j_2=0$, we only use the strong continuity. For $j_1>0$ and $j_2=0$, we use the cancellation with $B_{Q_1}$ : by the classical arguments, we have that for $x\in S_{j_1}(Q_1) \cap Q_2$ and $g=B_{Q_2}(f_2)$
$$\left|T(B_{Q_1}(f_1),g)(x)\right|\leq \frac{1}{\mu(Q_1)}\int_{Q_1} \int_{Q_1} \int_X \left| k(x,y,z)-k(x,h,z)\right| |f_1(y)||g(z)| d\mu(z) d\mu(y) d\mu(h).$$
By using the estimate about the kernel, we obtain~:
$$\left|T(B_{Q_1}(f_1),g)(x)\right|\leq \left(\int_X \frac{\mu(Q_1)^{1+\epsilon}}{\mu(2^{j_1}Q_1) + \mu(B(x,d(z,x)))^{2+\epsilon}} d\mu(z) \right) \|f_1\|_\infty \|f_2\|_\infty.$$
We can compute the integral over the variable $z\in X$ and then we get 
$$\left|T(B_{Q_1}(f_1),g)(x)\right|\leq  \frac{\mu(Q_1)^{\epsilon}}{\mu(2^{j_1}Q_1) )^{1+\epsilon}} \|f_1\|_\infty \|f_2\|_\infty \lesssim 2^{-j\delta \epsilon}\frac{\mu(Q_1)}{\mu(2^{j_1}Q_1)} \|f_1\|_\infty \|f_2\|_\infty.$$
This inequality implies what we want with $\gamma_j\lesssim 2^{-j\delta\epsilon}$. When the two index $j_1,j_2$ are non vanishing, we use the cancellation for the smallest ball and then we obtain by the same arguments the desired inequality with $\gamma_j\lesssim 2^{-j\delta \epsilon/2}$.
So the assumptions of Proposition \ref{prop:bilineaire} are satisfied. Therefore we can apply Theorem \ref{thm:bilineaire1} and obtain the continuity from $H^1_{F,ato} \times H^1_{F,ato}$ into $L^{1/2}$. This result was already obtained by L. Grafakos and N. Kalton in \cite{GK} for $X=\R^n$. In \cite{GT}, L. Grafakos and R. Torres have obtained the weak type estimate $L^1 \times L^1$ into $L^{1/2,\infty}$ for the euclidean space.

\mb Now, we add the dual assumption of (\ref{cz1}) and (\ref{cz2}) : 
\be{cz3} \left|k(x,y,z)-k(h,y,z)\right| \lesssim \left( \frac{d(x,h)}{ d(x,y)+d(z,x)} \right)^\epsilon \frac{1}{\mu(B(x,d(y,x))) + \mu(B(x,d(z,x)))^2} \ee for $d(x,h)\leq d(x,y)/2+d(x,z)/2$.\\
Then this class of bilinear Calder\'on-Zygmund operators (satisfying (\ref{cz1}),(\ref{cz2}) and (\ref{cz3})) is stable by duality. By using Theorem \ref{thm:final} and combining it with the continuity from $H^1_{F,ato} \times H^1_{F,ato}$ to $L^{1/2}$ we deduce that for all exponents $\T,\T_1,\T_2$ such that $1<\T_1,\T_2\leq \infty$ and
$$0<\frac{1}{\T}=\frac{1}{\T_1}+\frac{1}{\T_2}<2$$
the operator $T$ can be continuously extended from $L^{\T_1} \times L^{\T_2}$ into $L^{\T}$.
So for these bilinear Calder\'on-Zygmund operators (under the three assumptions (\ref{cz1}), (\ref{cz2}) and (\ref{cz3})), one strong continuity implies all these continuities. 

\gb In addition in this case, using Proposition \ref{prop:bilineaire2}, we obtain the continuity of $T$ from $H^1_{F,ato} \times  L^{\T_1}$ into $L^{\T}$ for $1<\T_1\leq \infty$ satisfying
$$\frac{1}{\T}=\frac{1}{\T_1}+1.$$
Applying this to $T^{*1}$ with $\T_1=\infty$, we obtain that $T$ admits a continuous extension from $L^\infty_c \times L^\infty_c$ into $BMO$. Here we note $L^\infty_c$ for the set of compactly supported and bounded functions.

\mb All these results are already proved in \cite{GT} in the euclidean case and the authors have described a ``T(1)-bilinear Theorem'' for this kind of operators to obtain a criterion to get the important first strong continuity. 

\subsubsection{Continuities in weighted Lebesgue spaces.}

\gb Now we are interesting to weighted inequalities. In \cite{GT2}, L. Grafakos and R. Torres have studied some weighted estimates for such bilinear Calder\'on-Zygmund operators in the euclidean case. They succeed to obtain results to prove boundedness from $L^{\T_1}_{\omega_1} \times L^{\T_2}_{\omega_2}$ into $L^{\T}_{\omega}$ when the three weights are equal : $\omega_1=\omega_2=\omega$ and belong to a certain class using a pointwise Cotlar's inequality. Then in \cite{GM2}, L. Grafakos and J.M. Martell have described a multiple weight extrapolation theory and have shown that a bilinear Calder\'on-Zygmund operator is continuous from $L^{\T_1}_{\omega_1} \times L^{\T_2}_{\omega_2}$ into $L^{\T}_{\omega}$ for
$\omega_1\in {\mathbb A}_{\T_1}$, $\omega_2\in {\mathbb A}_{\T_2}$ and $\omega:=\omega_1^{\T/\T_1}\omega_2^{\T/\T_2}$.
These results are based on the pointwise estimates (\ref{cz1}) and (\ref{cz2}). If we just have local estimates on the kernel, we could not have Cotlar's inequality. So it is interesting to describe what our weighted results  (Subsection \ref{subsectionpoids}) gives for bilinear operators satisfying only similar local estimates.

\gb We could probably improve our weighted results with ideas of the weighted multilinear extrapolation results of \cite{GM2}.
Let us give an example.

\mb
Let $T$ be a bilinear Calder\'on-Zygmund operator satisfying one strong continuity. From the previous discussion and according to Remark \ref{remim}, we can deduce that $U_g:=T(.,g)$ satisfies Theorem \ref{theo4} for $g\in L^\infty$. By using Remark \ref{remi}, we can also deduce that for a weight $\omega_1\in {\mathbb A}_{p_1}$ the operator $U_g$ admits a continuous extension from $L^{p_1}_{\omega_1}$ into $L^{p_1}_{\omega_1}$ for an exponent $1<p_1<\infty$.
So we get that $T$ is continuous from $L^{p_1}_{\omega_1} \times L^\infty$ into $L^{p_1}_{\omega_1}$. By symmetry for a weight $\omega_2 \in {\mathbb A}_{p_2}$ with $1<p_2<\infty$, we obtain that $T$ is continuous from $L^\infty \times L^{p_2}_{\omega_2}$ into $L^{p_2}_{\omega_2}$. Now by bilinear interpolation, we obtain that $T$ is bounded from 
$L^{\T_1}_{\omega_1} \times L^{\T_2}_{\omega_2}$ into $L^{\T}_\omega$ for exponents $1<\T,\T_1,\T_2$ such that there exists $\theta\in]0,1[$ with
$$ \frac{1}{\T_1}= \frac{\theta}{p_1}, \qquad \frac{1}{\T_2}= \frac{1-\theta}{p_2}, \qquad \frac{1}{\T}= \frac{1}{\T_1}+\frac{1}{\T_2},$$
and with the weight $\omega=\omega_1^{p_1\theta/\T} \omega_2^{(1-\theta)p_2/\T}$. The real interpolation theory of weighted Lebesgue spaces is described in the book \cite{BL} at Sections 5.4 and 5.5. However this result is weaker as the one described in \cite{GM2}. Using their main weighted extrapolation result, from the two continuities $L^{p_1}_{\omega_1} \times L^\infty$ into $L^{p_1}_{\omega_1}$ and from $L^\infty \times L^{p_2}_{\omega_2}$ into $L^{p_2}_{\omega_2}$ we can regain their results described in \cite{GM2}. So it seems to be useful to make a mixture of our weighted results and the weighted extrapolation theory to obtain the strongest results.

\subsection{The generalized bilinear Calder\'on-Zygmund operators.}

\mb In \cite{DGY}, X.T. Duong, L. Grafakos and L. Yan have generalized the previous example in considering the same kind of bilinear operators associated to other oscillations.
We explain in this subsection, how we can regain their results.

\mb On $\R^n$, they choose operators $A_t$, which are given by their kernels $a_t$ satisfying
$$\left|a_t(x,y)\right|\leq t^{-n/s}h\left(\frac{|x-y|}{t^{1/s}}\right),$$
where $s>0$ is a fixed parameter and $h$ is a positive bounded decreasing function satisfying
$$ \lim_{r\to \infty} r^{sn+\eta} h(r^s) =0$$
for a parameter $\eta>0$. So the operator $A_t$ is uniformly $L^\infty$-bounded and uniformly $L^1$-bounded.

\mb According to these notations, we define for $Q=B(x,r)$ a ball the operator $A_Q:=A_{r^{1/s}}$, and then $B_Q=Id-A_Q$. It is obvious to check that the assumptions on the kernel $a_t$ give us that our maximal operator $M_{\infty}$ (defined by (\ref{opeM})) is bounded by the maximal function $M_{HL,1}$. So with the exponent $\B=\infty$, we can define our Hardy space $H^1_{ato}$, which is $L^1-L^\infty$ regularizing. 

\mb In \cite{DGY}, the authors used two Assumptions~: Assumption 2.1 and Assumption 2.2, we recall. The bilinear operator $(f,g) \rightarrow T(B_Q(f),g)$ has a bilinear kernel $k_Q^1(x,y,z)$ satisfying the following estimate
$$ \left|k_Q^1(x,y,z)\right| \lesssim \frac{r_Q^\epsilon}{(|x-y|+|x-z|)^{2n+\epsilon}} + \frac{{\bf 1}_{[-r_Q,r_Q]}(y-z)}{(|x-y|+|x-z|)^{2n}}$$
for $|x-y|\geq r_Q$ and similarly for the operator $(f,g) \rightarrow T(f,B_Q(g))$.
These assumptions permit to the operator $T$ to satisfy (\ref{bi1}) of our Proposition \ref{prop:bilineaire} when one the index $j_1,j_2$ is non vanishing. The proof is similar to the one of the previous example, here there is an extra term to be studied but similar arguments permit to obtain the desired estimates (we let the details to the reader). So as for the previous example, if we already have one strong continuity for this operator $T$, then we have the continuity in $H^1_{F,ato} \times H^1_{F,ato}$ into $L^{1/2}$, which can be compared to the weak type estimate $L^1\times L^1 \to L^{1/2,\infty}$ obtained in \cite{DGY}. Then by the assumptions 3.1 and 3.2 of \cite{DGY} (which correspond to the previous assumptions for the operator $T$ and its two adjoints), we construct a class of bilinear operators which is stable by duality and so we can apply our Theorem \ref{thm:final} to obtain a new proof of Theorem 3.1 of \cite{DGY}.

\mb Using the same arguments as for the previous example, we can deduce that such operators admit a continuous extension from $L^\infty_c \times L^\infty_c$ into the space $BMO_A$, which is defined by the norm
$$\|f\|_{BMO_A}:=\sup_{Q} \left(\frac{1}{\mu(Q)}\int_Q \left|f-A_Q^*(f) \right| d\mu\right).$$
To prove this claim, we use Proposition 9.2 of \cite{BJ}, which ``characterizes'' the dual space $(H^1_{F,ato})^*\cap L^2$ by this $BMO$ space. This end-point estimate seems to be new compared to results in \cite{DGY}. 

\gb So our results permit us to obtain a new proof of the main theorem of \cite{DGY}. We have proved our results with the most abstract framework and assumptions and so we can generalize this example with restriction for exponents, more general operators $A_t$ and work on a space of homogeneous type. 
\mb In addition as for the previous example, we can develop a multiple weight theory for these operators. In \cite{DGY}, the authors show how multilinear Calder\'on commutators could be thought as a particular case of these generalized bilinear Calder\'on-Zygmund operators.

\subsection{Applications to quadratic functionals.}

In this subsection, we would like to describe how we can use this bilinear theory to study quadratic linear functionals.
For example let $X=\R^n$ and $A$ be an $n\times n$ matrix-valued function satisfying the ellipticity condition~: there exist two constants $\Lambda\geq \lambda>0$ such that
$$ \forall \xi,\zeta\in \C^n, \qquad \lambda|\xi|^2 \leq Re \left( A\xi \cdot \overline{\xi} \right) \quad  \textrm{and} \quad |A\xi \cdot \overline{\zeta} | \leq \Lambda|\xi||\zeta|.$$ 
We define the second order divergence form operator
$$ L(f):= -\textrm{div} (A \nabla f),$$
and then we compute the quadratic linear functional~:
 $$ S_L(f)(x):=\left(\int_{0}^\infty \left|(tL)^{1/2} e^{tL} f(y) \right|^2 \frac{dt}{t} \right)^{1/2}. $$ 
We define the limit exponent $p_-$ as
$$ p_-:= \inf \left\{p\geq 1, \ \sup_{t>0} \|e^{tL}\|_{L^p\to L^p} <\infty \right\} .$$
In \cite{A}, P. Auscher have proved that for $p\in (p_-,2]$, the sublinear operator $S_L$ is $L^p$ bounded. We describe how we can regain this result. We bilinearize the square function with the following bilinear operator~:
$$BS_L(f,g)(x):= \int_{0}^\infty (tL)^{1/2} e^{tL} f(y) (tL)^{1/2} e^{tL} g(y) \frac{dt}{t}.$$
We have the direct equivalence for $1<p<\infty$ : $S_L$ is $L^p$ bounded if and only if $BS_L$ is bounded on $L^{p} \times L^p$ in $L^{p/2}$.
By the $L^2$ functional calculus, we start from the $L^2$ boundedness of $S_L$ so $BS_L$ is bounded from $L^2 \times L^2$ into $L^1$. Now we will use an adapted Hardy space to obtain other continuities~: for every ball $Q$, we choose our oscillation operator
$$ B_Q(f) := (r_Q^2L)^{M}e^{-r_Q^2L}(f) \qquad \textrm{or} \qquad   B_Q(f) := \left(Id - (Id+r_Q^2L)^{-1} \right)^{M}(f).$$
Then we can construct the atomic Hardy space $H^1_{F,ato}$.
Using the $L^2$ off-diagonal decays of the semigroup (see Section 2.3 of \cite{A}), it is quite ``classical'' to obtain the assumptions of our Proposition \ref{prop:bilineaire} (the arguments are very similar to those of Theorem 6.1 in \cite{A} (step3) and those in \cite{HM}). Then we get the boundedness of $BS_L$ from  $H^1_{F,ato} \times H^1_{F,ato}$ into $L^{1/2}$.
In addition for all $p_0\in (p_-,2)$ we have $L^{p_0}-L^2$ off diagonal estimates for the semi group $e^{tL}$. Using Remark \ref{remms}, we know that our Hardy space $H^1_{F,ato}$ is $L^{p_0}-L^2$ regularizing. Applying Theorem \ref{thm:bilineaire1}, we also deduce that for all exponent $p\in (p_-,2]$ our bilinear operator $BS_L$ admits a continuous extension from $L^p \times L^p$ into $L^{p/2}$, which is equivalent to the $L^p$ boundedness of $S_L$.

\gb We have let the details to the interested reader. Here we want just explain how use the bilinear theory to study quadratic functionals. It is interesting to note that this point of view permit us to obtain the desired result without resorting $l^2$ valued arguments.

\subsection{The bilinear Marcinkiewicz multipliers.}

\mb In \cite{GK2}, L. Grafakos and N. Kalton have studied bilinear Marcinkiewicz multipliers on $\R^n$. A bilinear operator $T$ is a bilinear Marcinkiewicz multiplier if it is associated to a symbol $\sigma$ by
$$T(f,g)(x):= \int_{\R^{2n}} e^{ix(\alpha+\beta)} \widehat{f}(\alpha) \widehat{g}(\beta) \sigma(\alpha,\beta) d\alpha d\beta$$ 
satisfying
$$\forall a,b\in \N^{n}, \qquad \left| \partial_\alpha^a \partial_\beta ^b \sigma(\alpha,\beta) \right| \leq |\alpha|^{-|a|} |\beta|^{-|b|}.$$
These operators are a little more singular than bilinear Calder\'on-Zygmund operators. Such an operator has a bilinear kernel $k(x-y,x-z)$ satisfying
\be{marcin} \forall a,b\in \N^{n}, \qquad \left| \partial_y^a \partial_z ^b k(x-y,x-z) \right| \leq |x-y|^{-n-|a|} |x-z|^{-n-|b|}. \ee

\gb We want in this subsection to explain what results can we obtain for these bilinear operators. 

\gb From (\ref{marcin}), a bilinear Marcinkiewicz multiplier is almost a bilinear Calder\'on-Zygmund operator. So as for the first example, let us take our operators $B_Q$ equal to the exact oscillation operator
 $$B_Q(f):= f(x) - \left(\frac{1}{\mu(Q)}\int_Q f d\mu \right) {\bf 1}_{Q}.$$

\mb We would like to apply Proposition \ref{prop:bilineaire}. Let us take its notations. When the two index $j_1,j_2$ are non vanishing, we can use cancellation on the two balls and as for the classical bilinear Calder\'on-Zygmund operators, we can find very fast decay and so we can choose for the coefficients $\gamma_j:=2^{-jn/2}$. When the two indexes are equal to $0$, we will use as previously a strong continuity (assumed on the bilinear operator) to obtain the desired inequality. The main difficulty is when one of the two index is equal to $0$.

\mb In \cite{GK2}, the authors obtained an equivalent condition over the symbol $\sigma$ to a strong continuity for the operator $T$. It is interesting to note that this condition is independent on the different exponents. So they proved that if such an operator is continuous from $L^{\T_1}\times L^{\T_2}$ into $L^{\T}$ with $1<\T_1,\T_2<\infty$ and
$$ \frac{1}{\T}=\frac{1}{\T_1}+\frac{1}{\T_2}$$
then it is continuous for all the exponents satisfying the same properties.
Due to the similarity, we probably can show the assumptions of Proposition \ref{prop:bilineaire} when one of the two index $j_1,j_2$ is equal to $0$, assuming one strong continuity. But today this fact is not clear for us.

\mb We finish to underline the improvement of these ideas, compared with the ideas based on a Calder\'on-Zygmund decomposition (used in \cite{DGY,GT,GT2}). The use of an appropriate Hardy space and the bilinear interpolation theory permit us to reduce the problem to the linear theory of Hardy spaces. In addition, it permits to study kindly the two arguments of bilinear operators and so we hope to deduce a multiple weight theory for these kind of operators. It will be interesting to combine our results with those of the weighted extrapolation theory. In addition, we have given in Proposition \ref{prop:bilineaire}, a criterion for a bilinear operator to act on a Hardy space $H^1_{F,ato}$. It will be interesting to obtain a weaker condition, for example we have just seen that for a Marcinkiewicz multiplier this condition is not obvious.

\bibliographystyle{plain}
\bibliography{hardy}

\end{document}